\documentclass[final,1p,times]{elsarticle}
\usepackage{graphicx}
\usepackage{amssymb,amsmath,bm}  

\def\E{{\bf E}}
\def\H{{\bf H}}
\newcommand{\bi}{\begin{itemize}}
\newcommand{\ei}{\end{itemize}}
\newcommand{\ben}{\begin{enumerate}}
\newcommand{\een}{\end{enumerate}}
\newcommand{\be}{\begin{equation}}
\newcommand{\ee}{\end{equation}}
\newcommand{\bea}{\begin{eqnarray}}
\newcommand{\eea}{\end{eqnarray}}
\newcommand{\bc}{\begin{center}}
\newcommand{\bfi}{\begin{figure}}
\newcommand{\efi}{\end{figure}}
\newcommand{\ca}[2]{\caption{#1 \label{#2}}}
\newcommand{\ig}[2]{\includegraphics[#1]{#2}}
\newcommand{\ec}{\end{center}}

\newcommand{\tbox}[1]{{\mbox{\tiny #1}}}
\newcommand{\mbf}[1]{{\bm #1}}           
\newcommand{\pO}{{\partial\Omega}}
\newcommand{\pU}{{\partial U}}
\newcommand{\eps}{\varepsilon}
\newcommand{\emach}{\eps_\tbox{mach}}
\newcommand{\bp}{{\bf Proof:\ }}   
\newcommand{\ep}{\hfill $\square$ \vspace{1ex} \\} 
\newcommand{\xx}{\textbf{x}}             
\newcommand{\yy}{\textbf{y}}
\newcommand{\zz}{\textbf{z}}
\newcommand{\kk}{\textbf{k}}
\newcommand{\dd}{\textbf{d}}
\newcommand{\qq}{\textbf{q}}
\newcommand{\rr}{\textbf{r}}
\newcommand{\eb}{\textbf{e}_1}                 
\newcommand{\el}{\textbf{e}_2}
\newcommand{\gqp}{G_{\tbox{\rm QP}}}              
\newcommand{\gqpr}{G_{\tbox{\rm QP}}^r}           
\newcommand{\gqprr}{\tilde{G}_{\tbox{\rm QP}}^r}  
\newcommand{\aqp}{A_{\tbox{\rm QP}}}
\newcommand{\aqpm}{\tilde{A}_{\tbox{\rm QP}}}    
\newcommand{\om}{\omega}
\newcommand{\al}{\alpha}
\newcommand{\bt}{\beta}
\newcommand{\vt}[2]{\left[\begin{matrix}#1\\#2\end{matrix}\right]}
\newcommand{\vf}[4]{\left[\begin{matrix}#1\\#2\\#3\\#4\end{matrix}\right]}
\newcommand{\mt}[4]{\left[\begin{matrix}#1&#2\\#3&#4\end{matrix}\right]}
\newcommand{\un}{{u_n}}
\newcommand{\RR}{\mathbb{R}^2}

\DeclareMathOperator{\Null}{Null}
\DeclareMathOperator{\dist}{dist}
\newtheorem{thm}{Theorem}
\newtheorem{lem}[thm]{Lemma}
\newtheorem{dfn}[thm]{Definition}
\newtheorem{rmk}[thm]{Remark}

\newtheorem{cnj}[thm]{Conjecture}
\newcommand{\matlab}{MATLAB}       

\journal{J. Comput. Phys.}

\begin{document} 
\begin{frontmatter}

\title{A new integral representation for quasiperiodic fields
and its application to two-dimensional band structure calculations}

\author[AHB]{Alex Barnett\corref{c:AHB}}
\ead{ahb@math.dartmouth.edu}
\ead[url]{http://www.math.dartmouth.edu/$\sim$ahb}
\author[LG]{Leslie Greengard}
\ead{greengard@cims.nyu.edu}
\ead[url]{http://math.nyu.edu/faculty/greengar}
\cortext[c:AHB]{Corresponding author. tel:+1-603-646-3178. fax:+1-603-646-1312}
\address[AHB]{Department of Mathematics, 
Dartmouth College, Hanover, NH, 03755, USA}
\address[LG]{Courant Institute, New York University, 251 Mercer St,
NY, 10012, USA}

\begin{abstract}
In this paper, we consider 
band-structure calculations governed by the Helmholtz or Maxwell 
equations in piecewise homogeneous periodic materials.
Methods based on boundary integral equations are natural in this context,
since they discretize the interface alone and
can achieve high order accuracy in complicated geometries.
In order to handle the {\em quasi-periodic} conditions which are
imposed on the unit cell, the free-space Green's function is typically
replaced by its quasi-periodic cousin.
Unfortunately, the quasi-periodic Green's function diverges
for families of parameter values that correspond to resonances
of the empty unit cell.
Here, we bypass this problem by means of a new integral
representation that relies on the free-space Green's function alone,
adding auxiliary layer potentials on the boundary of the unit cell itself.
An important aspect of our method is that by carefully including
a few neighboring images, the densities may be kept smooth
and convergence rapid.
%
This framework
results in an integral equation  of the second kind, avoids spurious
resonances, and achieves spectral accuracy.
Because of our image structure, inclusions which intersect 
the unit cell walls may be handled easily and automatically.
Our approach is compatible with fast-multipole acceleration,
generalizes easily to three dimensions, and
avoids the complication of divergent lattice sums.
\end{abstract}

\end{frontmatter}

\section{Introduction}
\label{s:i}

A number of problems in wave propagation require the calculation of 
{\em quasi-periodic}
solutions to the governing partial differential equation in the frequency
domain. 
For concreteness, let us consider the two-dimensional  (locally isotropic)
Maxwell equations in what is called TM-polarization \cite{jackson,jobook}.
In this case, the Maxwell equations reduce to a scalar Helmholtz equation
 \be
\Delta u(x,y) +\omega^2 \epsilon \mu \, u(x,y) = 0,
\label{e:helm}
\ee
where 
$\epsilon$ and $\mu$ are the permittivity and permeability of the medium, respectively,
and we have assumed a time dependence of $e^{-i \omega t}$ at frequency $\omega > 0$.
Given a solution $u$ to \eqref{e:helm}, it is straightforward to verify that the corresponding
electric and magnetic fields $\E,\H$ of the form
 \begin{eqnarray*}
 \E(x,y,z) &=& \E(x,y) = \,
(0, 0, u(x,y))\\
 \H(x,y,z) &=& \H(x,y) = \, \frac{1}{i \omega \mu}
(u_y(x,y), -u_x(x,y), 0)
\end{eqnarray*}
satisfy the full system
\begin{eqnarray*}
\nabla \times \E &=& i \omega \mu \H \\
\nabla \times \H &=& -i \omega \epsilon \E \, .
\end{eqnarray*}
We are particularly concerned with doubly periodic materials whose
refractive index $n = \sqrt{\epsilon \mu}$ is piecewise constant
(Fig. \ref{f:g}). Such structures are typical in solid state physics,
and are of particular interest at present because of the potential
utility of photonic crystals, where the obstacles are dielectric
inclusions with a periodicity on the scale of the wavelength of light
\cite{jobook}. Photonic crystals allow for the control of optical wave
propagation in ways impossible in homogeneous media, and are finding a
growing range of exciting applications to optical devices, filters
\cite{channeldrop}, sensors, negative-index and meta-materials
\cite{metareview}, and solar cells \cite{PCsolar}.

\bfi [t]
a)\quad\raisebox{-1.8in}{\ig{width=0.45\linewidth}{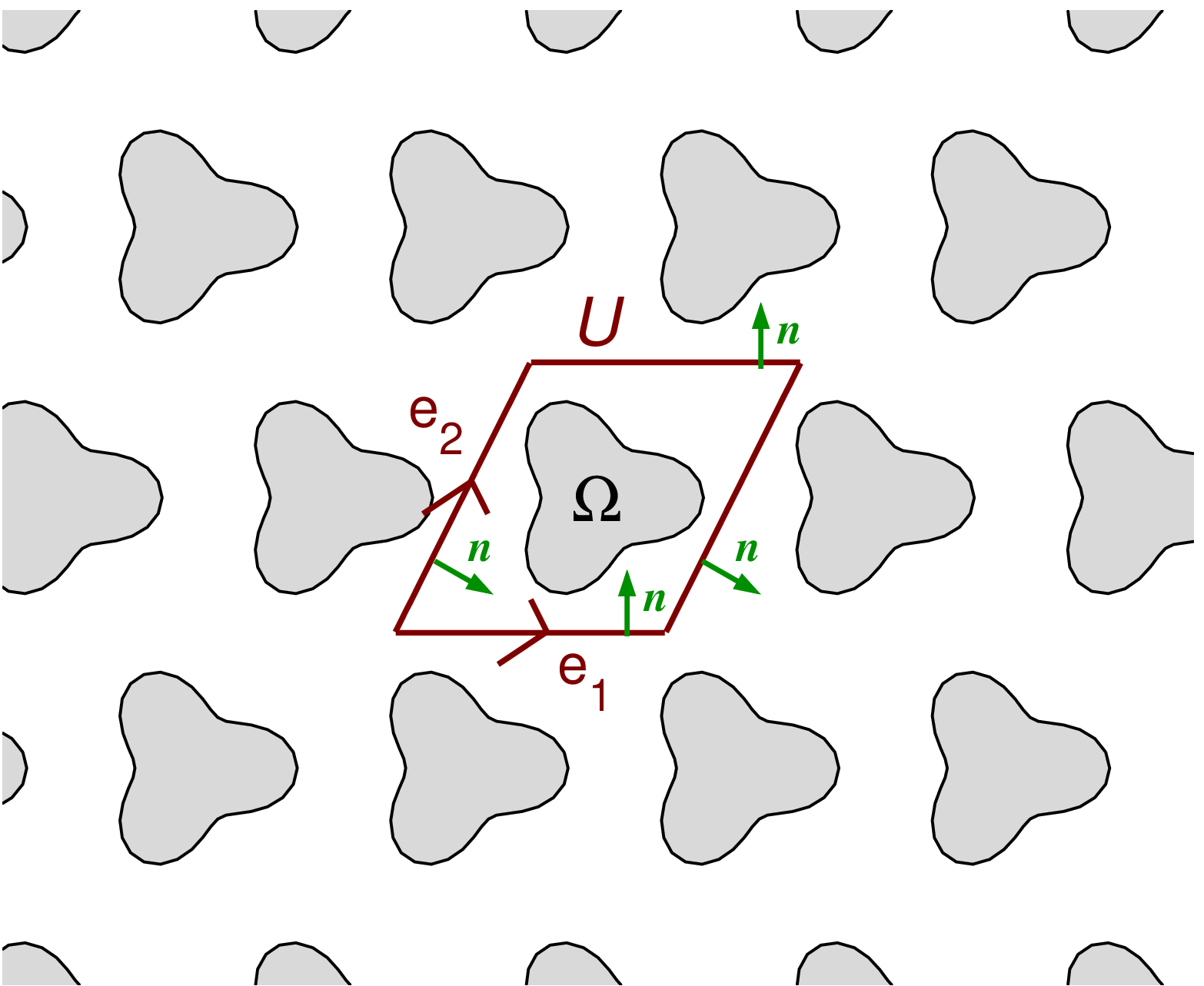}}
\quad
b)\quad\raisebox{-1.8in}{\ig{width=0.45\linewidth}{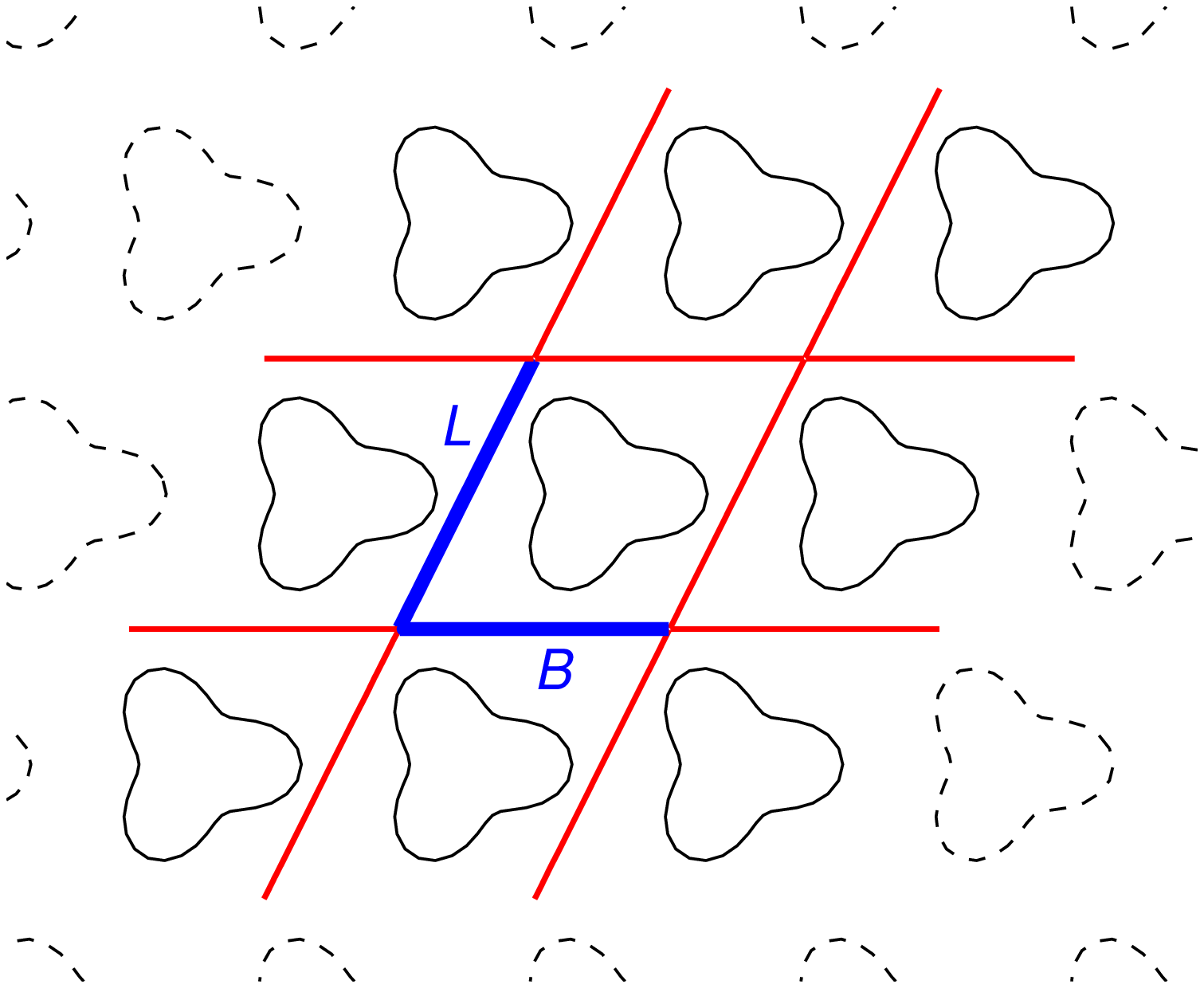}}
\ca{a) Problem geometry: an infinite dielectric crystal, in the case where
the inclusion $\Omega$ lies within a parallelogram unit cell $U$.
The (shaded) set of all inclusions in the lattice, denoted by $\Omega_\Lambda$
in the text, has refractive index $n$, while the white region has index 1.
b) Sketch of our quasi-periodizing scheme: we make use of
layer potentials on the left (L) and bottom (B) walls,
extended to the additional
segments shown, which form a skewed `tic-tac-toe' board,
as well as the near neighbor images of $\Omega$, outlined in solid lines.
}{f:g}
\efi

We assume that the crystal consists of a periodic array of obstacles
($\Omega_\Lambda$) with refractive index $n \neq 1$, embedded in a
background material with refractive index $n=1$ (denoted by
$\RR\setminus\overline{\Omega_\Lambda}$).  We then rewrite
(\ref{e:helm}) as a system of Helmholtz equations
\bea
(\Delta+n^2\omega^2)u &=& 0 \qquad \mbox{in } \Omega_\Lambda
\label{e:pdei}
\\
(\Delta+\omega^2)u &=& 0 \qquad \mbox{in } \RR\setminus\overline{\Omega_\Lambda}
\label{e:pdee}
\eea
The expression $\overline{\Omega_\Lambda}$, above, is used to denote
the closure of the domain $\Omega_\Lambda$ (the union of the domain
and its boundary $\partial\Omega_\Lambda$).  In this formulation, we
must also specify conditions at the material interfaces.  These are
derived from the required continuity of the tangential components of
the electric and magnetic fields across $\partial\Omega_\Lambda$
\cite{jackson,jobook}, yielding
\be
u, \un \mbox{ continuous across } \partial\Omega_\Lambda
\label{e:cont}
\ee
where $\un = \partial u/\partial n$ is the outward-pointing normal 
derivative. 

The essential feature of doubly periodic microstructures in 2D (or
triply periodic microstructures in 3D) is that, at each frequency,
there may exist traveling wave solutions (Bloch waves) propagating in
some direction defined by a vector ${\bf k}$.

\begin{dfn}
{\em Bloch waves} are nontrivial solutions to 
\eqref{e:pdei}--\eqref{e:cont} that are quasiperiodic, in the sense that
\be
 u(\xx) = e^{i\kk\cdot\xx}\tilde{u}(\xx) \, , 
\label{e:qpbloch}
\ee
where $\tilde{u}$ is periodic with the lattice period and $\kk=(k_x,k_y)$ is
real-valued. $\kk$ is referred to as the {\em Bloch wavevector}.
\end{dfn}

Bloch waves characterize the bulk optical properties at frequency
$\om$; they are analogous to plane waves for free space.  If such
waves are absent for all directions $\kk$ for a given $\omega$, then
the material is said to have a {\em band-gap} \cite{yablonovich}). The
{\em size} of a band-gap is the length of the frequency interval
$[\om_1,\om_2]$ in which Bloch waves are absent.  Crystal structures
with a large band-gap are `optical insulators' in which defects may be
used as guides \cite{jobook}, with the potential for enabling
high-speed integrated optical computing and signal processing.

\begin{dfn}
The {\em band-structure} of a given crystal geometry is the set of 
parameter pairs $( \omega,\kk )$ for which nontrivial Bloch waves exist. 
 \end{dfn}

 The numerical prediction of band structure is a computationally
 challenging task, yet essential to the design and optimization of
 practical devices. It requires characterizing the nontrivial
 solutions to a homogeneous system of partial differential equations
 \eqref{e:pdei}, \eqref{e:pdee} subject to homogeneous interface and
 periodicity conditions \eqref{e:cont}, \eqref{e:qpbloch}
 in complicated geometry.  Solving this eigenvalue problem is the
 focus of our paper.

In the next section, we briefly review existing approaches, and in section
\ref{s:gqp},
we present and test a method that relies on the quasi-periodic Green's function.
We introduce our new mathematical formulation in section \ref{s:qplp}. 
Numerical results are presented in section \ref{s:res}, and we conclude in 
section \ref{s:conc} with some remarks about the potential for 
wider application of this approach.

\section{Existing approaches}

In order to pose the band-structure problem as an eigenvalue problem
on the unit cell $U$ (see Fig.~\ref{f:g}), we will require some
additional notation.  The nonparallel vectors $\eb,\el\in\RR$ define a
{\em Bravais lattice} $\Lambda:=\{m\eb+n\el : m,n\in\mathbb{Z}\}$.
Given a smooth, simply connected inclusion
$\Omega\subset\mathbb{R}^2$, we may formally define the corresponding
dielectric crystal by $\Omega_\Lambda := \{\Omega+\dd:
\dd\in\Lambda\}$. As indicated above, we assume that $\Omega_\Lambda$
has refractive index $n \neq 1$, and that the background
$\RR\setminus\overline{\Omega_\Lambda}$ has refractive index 1.  For the moment,
we assume that $\overline{\Omega}\subset U$ as illustrated in
Fig.~\ref{f:g}.  We will discuss the case of $\Omega$ crossing
$\partial U$ in Section~\ref{s:wrap}.

The quasi-periodicity condition (\ref{e:qpbloch}) can be rewritten as
a set of boundary conditions on the unit cell $U$, coupling the
solution on the left ($L$) and right ($L+\eb$) walls, as well as on the bottom
($B$) and top ($B+\el$) walls. More precisely, if we define
\[a := \kk \cdot \eb, \quad \al:= e^{ia}, \qquad
b := \kk \cdot \el, \quad \bt:=e^{ib} \, ,  \]
then quasi-periodicity is written
\bea
u|_{L+\eb} &=& \alpha \, u|_L
\label{e:f}
\\
\un|_{L+\eb} &=& \alpha \, \un|_L 
\\
u|_{B+\el}&=&  \beta \,  u|_B
\\
\un|_{B+\el} &=& \beta \, \un|_B \, ,
\label{e:gp}
\eea
where the normals have the senses shown in Fig.~\ref{f:g}.

The homogeneous equations \eqref{e:pdei}-\eqref{e:cont},
\eqref{e:f}-\eqref{e:gp} define
a partial differential equation (PDE)
eigenvalue problem on the torus $U$. By convention, the
band structure or Bloch eigenvalues
are generally defined as the subset of the parameter space
$\{(\omega,a,b): \omega>0, -\pi\le a<\pi, -\pi\le b<\pi\}$
for which nontrivial solutions $u:U\to\mathbb{C}$ exist.
The earlier definition of band-structure, based on 
(\ref{e:qpbloch}), allows for arbitrary values of 
$\kk$. It is clear, however, that one only needs
to consider a single period of $\kk$'s projection onto
 $\eb,\el$, which we have denoted
by $a,b$, to characterize the entire
set $(\omega,\kk)$ of nontrivial Bloch waves. This domain
$\{(a,b): -\pi\le a<\pi, -\pi\le b<\pi\}$ is (essentially)
what is referred to as the {\em Brillouin zone}.

Because the PDE is elliptic and $U$ is compact, for each $\kk$
there is a discrete set of eigenvalues 
$\{\omega_j(\kk)\}_{j=1}^\infty$, counting multiplicity, accumulating only
at infinity. 
Each $\omega_j(\kk)$ is continuous in $\kk$,
so that the bands form {\em sheets}.

Popular numerical methods for band structure calculations are reviewed in
\cite{jobook}. Broadly speaking, they may be classified
as either time-domain or frequency domain schemes.
In the first case, an initial pulse is evolved via the full wave equation
(typically using a finite-difference or finite-element approximation). If the 
simulation is sufficiently long, Fourier transformation in the time variable 
then reveals the full band structure.
In the second case, 
the eigenvalue problem  \eqref{e:pdei}-\eqref{e:cont},
\eqref{e:f}-\eqref{e:gp} 
is discretized directly.
Such frequency domain schemes can be further categorized as
\begin{enumerate}
\item PDE-based methods, which involve discretizing the unit cell
  using finite difference or finite element methods
  \cite{axmann99,dobson99,dossou},
\item plane-wave methods which expand the function $\tilde{u}$ in
  (\ref{e:qpbloch}) as a Fourier series, and apply the partial
  differential operator in Fourier space \cite{jobook,johnson01},
\item semi-analytic multipole expansion methods which apply largely
  to cylindrical or spherical inclusions \cite{botten03,Pissoort},
\item methods which use a basis of particular solutions to the PDE at
  a given frequency $\omega$ and enforce both interface and boundary
  conditions as a linear system, such as the ``multiple multipole'' or
  ``transfer-matrix'' method \cite{hafner90,smajic}, and
\item boundary integral (boundary element) methods \cite{yuan08},
  which includes the method described here.
\end{enumerate}

For a fixed $\kk$, methods of type (1) and (2) result in large, sparse
generalized eigenvalue problems whose lowest few eigenvalues
approximate the first few bands $\omega_j(\kk)$. They have the
advantage that they couple easily to existing robust linear algebraic
techniques. PDE-based methods, however, require discretization of the
entire cell in a manner that accurately resolves the geometry of the
inclusion $\Omega$.  Plane-wave methods, which perform extremely well
when the index of refraction $n$ is smooth, have low order convergence
when $n$ is piecewise constant, as in the present setting. Both
require a large number of degrees of freedom.

Methods of type (3), (4) or (5), on the other hand, represent the
solution using specialized functions (solutions of the PDE) whose
dependence on $\om$ is nonlinear.  As a result, they can be much more
efficient and high-order accurate, dramatically reducing the number of
degrees of freedom required. Unfortunately, however, they result in a
nonlinear eigenvalue problem involving all the parameters $\om$, $a$
and $b$, and somewhat non-standard techniques are required to find
values of the parameters for which the system of equations is singular
\cite{spence}.
%

We are particularly interested in using boundary integral methods
(BIEs), since they easily handle jumps in the index in complicated
geometry, have a well understood mathematical foundation, and can
achieve rapid convergence, limited only by the order of accuracy of
the quadrature rules used.  High order accuracy is important, not
only because of the reduction in the size of the discretized problem,
but in carrying out subsequent tasks, such as sensitivity analyses
\cite{crutchfield} through the numerical approximation of derivatives,
and the computation of band slopes (group velocity), and band
curvatures (group dispersion).

There is surprisingly little historical literature on using BIE for
band structure calculations, although the last few years have begun to
see some activity in this direction (see, for example,
\cite{yuan08}). There is, however, an extensive literature on integral
equations for {\em scattering} from periodic structures, which we do
not seek to review here. For some recent work and additional
references, see \cite{CMS,otani08}.

\section{Integral equations based on the quasi-periodic Green's function}
\label{s:gqp}

An elegant approach to designing integral representations for
quasiperiodic fields involves the construction of the Green's function
that imposes the desired conditions (\ref{e:f})-(\ref{e:gp})
exactly. We first need some definitions \cite{CK83,m+f2}.
At wavenumber $\om>0$, the free space
Green's function for the Helmholtz equation, $G$ is defined by
$-(\Delta + \omega^2) G = \delta_{\mbf{0}}$ where $\delta_\mbf{0}$ is
the Dirac delta function centered at the origin. In 2D, this yields
\be
G(\xx) = G^{(\om)}(\xx) = \frac{i}{4}H_0^{(1)}(\om|\xx|),
\qquad\xx\in\RR\setminus\{\mbf{0}\},
\ee
where $H_0^{(1)}$ is the
outgoing Hankel function of order zero.  By formally summing over
images of the Green's function placed on the lattice $\Lambda$, with
correctly assigned phases, we get an explicit expression for the
quasi-periodic Greens function
\be
\gqp(\xx) = \sum_{\dd\in\Lambda} e^{i\kk\cdot\dd}G(\xx-\dd)
= \sum_{m,n\in\mathbb{Z}} \al^m \bt^n G(\xx-m\eb-n\el) \, .
\label{e:gqp}
\ee
We leave it to the reader to verify that $\gqp$ does, indeed, satisfy 
(\ref{e:f})-(\ref{e:gp}). One small caveat: the series in 
(\ref{e:gqp}) is conditionally convergent for real $\omega$. The 
physically meaningful limit is taken by assuming some dissipation 
$\omega = \omega + i \eps$ in the limit $\eps \rightarrow 0^+$
(see \cite{dienst01} for a more detailed discussion).
It will be useful to distinguish between the copy of the Green's function
sitting in the unit cell $U$ and the set of all other images. For this, we define the 
``regular'' part of the quasi-periodic Green's function by 
\be
\gqpr(\xx) =  \sum_{\substack{m,n\in\mathbb{Z} \\
(m,n) \neq (0,0)}} \al^m \bt^n G(\xx-m\eb-n\el) \, .
\label{e:gqpr}
\ee
This function is a smooth solution to the Helmholtz equation
within $U$ and clearly satisfies
\be
\gqp(\xx) =  G(\xx) + \gqpr (\xx)~.
\label{e:gqpsplit}
\ee

A spectral representation also exists
\cite{Beylkin08,dienst01}, built from
the plane-wave eigenfunctions of the quasi-periodic torus $U$:
\be
\gqp(\xx) = \frac{1}{\mbox{Vol}(U)}\sum_{\qq\in \Lambda^*}
\frac{e^{i(\kk+\qq)\cdot\xx}}{|\kk+\qq|^2 - \om^2}  \, .
\label{e:gqpspec}
\ee
Here, $\Lambda^\ast := \{m \rr_1 + n \rr_2: m,n \in \mathbb{Z}\}$
is the {\em reciprocal lattice} with vectors
$\rr_j$ defined by ${\textbf{e}}_i\cdot\rr_j = 2\pi\delta_{ij}$ for
$i,j = 1,2$.
From the denominators in \eqref{e:gqpspec}
it is clear that $\gqp$ may blow up for specific combinations of 
$\om$ and $\kk$. The quasiperiodic Green's function is, in fact,
well-defined if and 
only if  those parameters satisfy the following non-resonance condition.
\begin{dfn}[empty resonance] 
A parameter set $(\om,\kk)$, equivalently $(\om,a,b)$, is
{\em empty resonant} if $\om = |\kk + \qq|$
for some $\qq \in \Lambda^\ast$, otherwise it is {\em empty non-resonant}.
\end{dfn}
Our terminology comes from the fact that the blow-up in $\gqp$
is physically the resonance of the `empty' 
unit cell $U$, with refractive index 1 everywhere and
quasi-periodic boundary conditions.
That is, $\gqp$ is undefined if and only if 
$(\om,a,b)$ lies on the band structure of the empty unit cell.
The blow-up of the Green's function is less apparent 
from (\ref{e:gqp}), but is manifested in the divergence of the series, even
in the limit $\omega = \omega + i \eps$ with $\eps \rightarrow 0^+$.

It will be convenient sometimes to refer to a Green's function as a
function of two variables, with $G(\xx,\yy) := G(\xx-\yy)$, and
$\gqp(\xx,\yy) := \gqp(\xx-\yy)$.  Then, for each $\yy\in \RR$, the
function $\gqp(\cdot,\yy)$ is quasi-periodic.

We now represent solutions to the PDE eigenvalue problem 
\eqref{e:pdei}-\eqref{e:cont}, \eqref{e:f}-\eqref{e:gp}
by the layer potentials,
\be
u \;=\; \left\{\begin{array}{ll}
{\cal S}^{(n\om)}\sigma + {\cal D}^{(n\om)}\tau & \mbox{in } \Omega\\
{\cal S}^{(\om)}_\tbox{QP}\sigma + {\cal D}^{(\om)}_\tbox{QP}\tau & \mbox{in }
U\setminus\overline{\Omega}
\end{array}\right.
\label{e:qprep}
\ee
where the usual single and double layer densities \cite{CK83}
at any wavenumber $\om>0$ are defined
by
\bea
({\cal S}^{(\om)}\sigma)(\xx) &=& \int_\pO G^{(\om)}(\xx,\yy) \sigma(\yy) ds_\yy
\label{e:s}
\\
({\cal D}^{(\om)}\tau)(\xx) &=& \int_\pO
\frac{\partial G^{(\om)}}{\partial n_\yy}(\xx,\yy) \tau(\yy) ds_\yy
\label{e:d}
\eea
and their quasi-periodized versions are likewise
\bea
({\cal S}^{(\om)}_\tbox{QP}\sigma)(\xx) &=& \int_\pO \gqp^{(\om)}(\xx,\yy) \sigma(\yy) ds_\yy
\label{e:sqp}         
\\
({\cal D}^{(\om)}_\tbox{QP}\tau)(\xx) &=& \int_\pO
\frac{\partial \gqp^{(\om)}}{\partial n_\yy}(\xx,\yy) \tau(\yy) ds_\yy ~.
\label{e:dqp}
\eea
Here $ds$ is the usual arc length measure on $\pO$, and the
derivatives are with respect to the second variable in the outward
surface normal direction at $\yy$.  It is clear \cite{CK83}
that the above four fields satisfy the Helmholtz equation at
wavenumber $\om$
%
in both $\Omega$ and $U\setminus \overline{\Omega}$.
Note that we have chosen a {\em non-periodized} representation within
the inclusion $\Omega$ in \eqref{e:qprep}, which has some analytic
advantages (see Theorem \ref{t:aqp} and the last paragraph in the Appendix).

Since $u$ in \eqref{e:qprep} satisfies \eqref{e:pdei}, \eqref{e:pdee},
and \eqref{e:f}-\eqref{e:gp},
all that remains is to solve for densities $\sigma$, $\tau$
such that the matching conditions \eqref{e:cont} are satisfied,
which we now address.

Using superscripts $+$ and $-$ to denote limiting values
on $\pO$, approaching from the positive and negative normal side respectively,
we use the field \eqref{e:qprep} and the
standard jump relations for single and double layer potentials
\cite{CK83,guentherlee} to write
\be
\vt{u^+-u^-}{\un^+-\un^-} \;=\; \left(\,
\mt{I}{0}{0}{I} + \mt{D^{(\om)}_\tbox{QP}-D^{(n\om)}}
{S^{(n\om)}-S_\tbox{QP}^{(\om)}}
{T_\tbox{QP}^{(\om)}-T^{(n\om)}}
{D^{(n\om)\,\ast}-D^{(\om)\,\ast}_\tbox{QP}}
\,\right)
\vt{\tau}{-\sigma}
\; =: \; \aqp \eta
\label{e:aqp}
\ee
Here $I$ is the identity operator, while $S$ and $D$ are defined to be the limiting
boundary integral operators (maps from $C(\pO)\to C(\pO)$)
with the kernels ${\cal S}$ and ${\cal D}$ interpreted in the 
principal value sense. ($S$ is actually weakly singular so the limit is already well defined. 
A standard calculation \cite{CK83,guentherlee} shows that $D$ is weakly singular as well).
The hypersingular operator $T$ has the kernel
$\frac{\partial^2 G(\xx,\yy)}{\partial n_\xx \partial n_\yy}$ and is unbounded as a map from
$C(\pO)\to C(\pO)$.
In these definitions, as in \eqref{e:s}-\eqref{e:dqp},
it is implied that $G$ inherits the appropriate superscripts and subscripts
from $S$, $D$ and $T$.
Finally, $\ast$ indicates the adjoint.
The amounts by which the material matching conditions fail to be
satisfied, 
\be
m \; := \; \vt{u^+-u^-}{\un^+-\un^-} ~,
\label{e:m}
\ee
is a column vector of functions which we call the {\em mismatch}.
We summarize the linear system
\eqref{e:aqp} by $m = \aqp \eta$ where $\eta:=[\tau; -\sigma]$.
It is important to note that the {\em difference} of hypersingular kernels,
${T_\tbox{QP}^{(\om)}-T^{(n\om)}}$,
in \eqref{e:aqp} is only weakly singular \cite[Sec.~3.8]{CK83}.
This cancellation, achieved here by using the same pair of densities
inside as outside the inclusion, is well known \cite{rokh83}.
The result is that $\aqp$ is a compact perturbation of the identity and
 \eqref{e:aqp} is a Fredholm system of integral equations of the second kind.

In the above scheme, we might hope that if
it is possible to find nontrivial densities $\eta$
whose field $u$ gives zero mismatch $m$ for a
set of parameters $(\om,a,b)$, then that set is a Bloch
eigenvalue. Indeed (as with the case of
simpler domain eigenvalue problems \cite[Sec.~8]{mitrea})
we have a stronger result. 
\begin{thm} 
Let $(\om,a,b)$ be empty non-resonant. Then
$(\om,a,b)$ is a Bloch eigenvalue if and only if
$\Null\aqp\neq\{0\}$~.
\label{t:aqp}
\end{thm} 
The proof occupies Appendix~\ref{a:aqp}. 
This suggests the core of a numerical scheme:
at each of a sampling (e.g.\ a grid) of parameters $(\om,a,b)$, find the
lowest singular value $\sigma_\tbox{min}(\aqpm)$
of a matrix discretization $\aqpm$ of $\aqp$.
The band structure will then be found where $\sigma_\tbox{min}(\aqpm)$
is close to zero.

\subsection{Discretization of the integral operators} 
\label{s:discr}

Since the goal of this work is to explore periodization, we limit ourselves
to the simplest case of $\pO$ being smooth. The methods of this paper extend 
without much effort to other shapes, but the quadrature issues become more involved.
Recalling \eqref{e:gqpsplit},
note that the kernels in \eqref{e:aqp} are the sum of a component
due to $G$ which is weakly singular, plus the remainder due to $\gqpr$
which is smooth (analytic). We will make use of a Nystr\"{o}m discretization using
the spectral quadrature scheme of Kress \cite{kress91} for $G$ and 
the trapezoidal rule for $\gqpr$.


We first remind the reader of the periodic trapezoidal Nystr\"{o}m scheme
\cite{LIE}, in the context of a general second kind boundary
integral equation 
\[ \mu(\xx) + \int_\pO k(\xx,\yy) \mu(\yy) ds_\yy = f(\xx),
\qquad \xx\in\pO, \]
where
$\pO$ is parametrized by the $2\pi$-periodic analytic function
$\zz: [0,2\pi) \to \RR$. Changing variable gives
\[ \mu(s) + \int_0^{2\pi} K(s,t) \mu(t) dt = f(s),
\qquad s\in[0,2\pi), \]
where $K(s,t):=k(\zz(s),\zz(t))\,|\zz'(t)|$ and $\zz'=d\zz/dt$.
Choosing $N$ quadrature points $t_j = 2\pi j/N$ with equal weights $2\pi/N$
gives the $N$-by-$N$ linear system for the unknowns
$\mu^{(N)}_j$, which approximate the exact values $\mu(t_j)$,
as
\be
\mu^{(N)}_k + \frac{2\pi}{N}\sum_{j=1}^N K(t_k,t_j)\mu^{(N)}_j
\;=\; f(t_k), \qquad k=1,\dots, N~.
\label{e:nyst}
\ee
By Anselone's theory of collectively compact operators \cite{LIE},
the convergence of errors $\bigl|\mu^{(N)}_j-\mu(t_j)\bigr|$ inherits the
order of the quadrature scheme applied to the exact integrand $K(s,\cdot)\mu$,
which is analytic when $k$ and $f$ are.
\begin{rmk}
For analytic integrands,
the periodic trapezoidal rule has exponential convergence with
error $O(e^{-2 \gamma N})$ where $\gamma$ is
the smallest distance from the real axis of any singularity
in the analytic continuation of the integrand.
\cite[Thm.~12.6]{LIE}.
\label{r:strip}
\end{rmk}

The above discretization is used to populate the matrix entries
in \eqref{e:aqp} that are due to the smooth compoment $\gqpr$.
(We explain how to compute this kernel itself in Section~\ref{s:eval}.)

For non-smooth kernels, such as $G$, the rule 
(\ref{e:nyst}) must be replaced by a quadrature that correctly
accounts for the singularity in order to retain high order accuracy. 
There are a variety of such schemes, such as those of \cite{alpert,helsing,kapur}.
By fixing the order of accuracy, they allow for straightforward coupling to 
fast multipole acceleration \cite{fmm1,fmm2,CMS,otani08} by making local modifications of 
a simple underlying quadrature rule (such as the trapezoidal rule or a composite
Gaussian rule). In the present context, we ignore such considerations and use a
global rule due to Kress \cite{kress91} that achieves spectral accuracy in the 
logarithmically singular case.

The essential idea of Kress' scheme (after transformation of variables to the 
interval $[0,2\pi]$) is to split a logarithmically singular kernel $K(s,t)$ in the form
\be
K(s,t) = \log\left(4 \sin^2 \frac{s-t}{2}\right) K_1(s,t) + K_2(s,t)
\label{e:Kst}
\ee
with $K_1$ and $K_2$ periodic and analytic. $K_2$ is (again) handled with the trapezoidal
rule. For $K_1$, the Kussmaul-Martensen
quadrature rule is spectrally accurate:
\be
\int_0^{2\pi} \log\left(4 \sin^2 \frac{s-t}{2}\right) g(t) dt
\approx
\sum_{j=1}^{N} R^{(N)}_j(s) g(t_j)
\ee
with quadrature weights (deriving from the Fourier series of the log factor)
given by
\be
R^{(N)}_j(s) \;=\; -\sum_{m=1}^{N/2-1}\frac{2}{m}
\cos m(s-t_j)\;-\;\frac{2}{N}\cos \frac{N}{2}(s-t_j)~.
\ee
Thus, the matrix elements in discretizing \eqref{e:Kst} are
$K(t_k,t_j) = R^{(N)}_{|j-k|}(0)K_1(t_k,t_j) + K_2(t_k,t_j)$.
Finally, it is always the difference of two hypersingular operators $T$ that
appears in the integral equation 
\eqref{e:aqp}. This difference is only logarithmically singular, so that Kress' rule
can be used for every block of \eqref{e:aqp}. We refer the reader to \cite{kress91}
for further details.

In summary, a matrix discretization $\hat{A}_\tbox{QP}$ of $\aqp$
is formed by using the above quadrature rules for each of the 
 2-by-2 integral operator blocks in \eqref{e:aqp}. This matrix 
maps density values to field values.
However, in order to create a matrix whose singular values
approximate those of $\aqp$ we must instead normalize such that
$2N$-dimensional Euclidean 2-norms correctly
approximate $L^2(\pO)$-norms.
This is done by symmetrizing using quadrature weights to give our final matrix
\be
\tilde{A}_\tbox{QP} = W^{1/2}\hat{A}_\tbox{QP}W^{-1/2}
\label{e:aqpw}
\ee
where $W$ is diagonal with diagonal elements
$w_j=w_{j+N}= (2\pi/N)|\zz'(t_j)|$, for $j=1,\ldots,N$. 

The net result of the preceding discussion is that with the use of specialized quadratures
on smooth boundaries,  the singular values of $\tilde{A}_\tbox{QP}$
are spectrally accurate approximations to those of $\aqp$.
We demonstrate this convergence for a small trefoil-shaped inclusion
in Fig.~\ref{f:gqpconv}a; the convergence is spectral, until
the error is approximately machine precision times the matrix 2-norm.
The rate appears to be faster at a Bloch eigenvalue
(in this case on the fourth band)
than far from one.
Fig.~\ref{f:gqpconv}b shows that the minimum locates the
parameter $b$ to 14 digit accuracy for $N\ge70$.


\bfi
\hspace{-5ex}
\mbox{%
a)\raisebox{-2.4in}{\hspace{-2ex}\ig{height=2.3in}{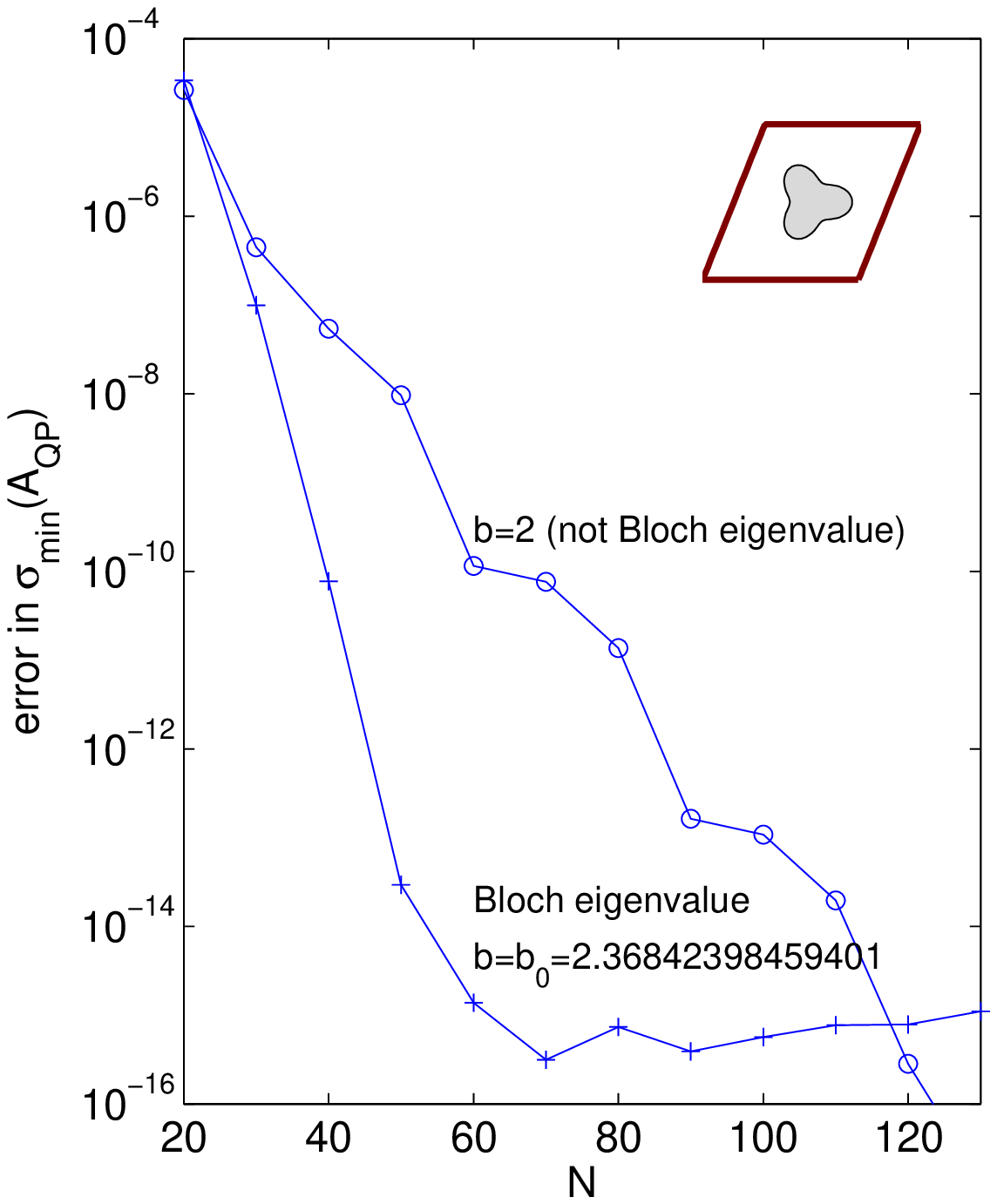}}
\;
b)\raisebox{-2.4in}{\hspace{-2ex}\ig{height=2.3in}{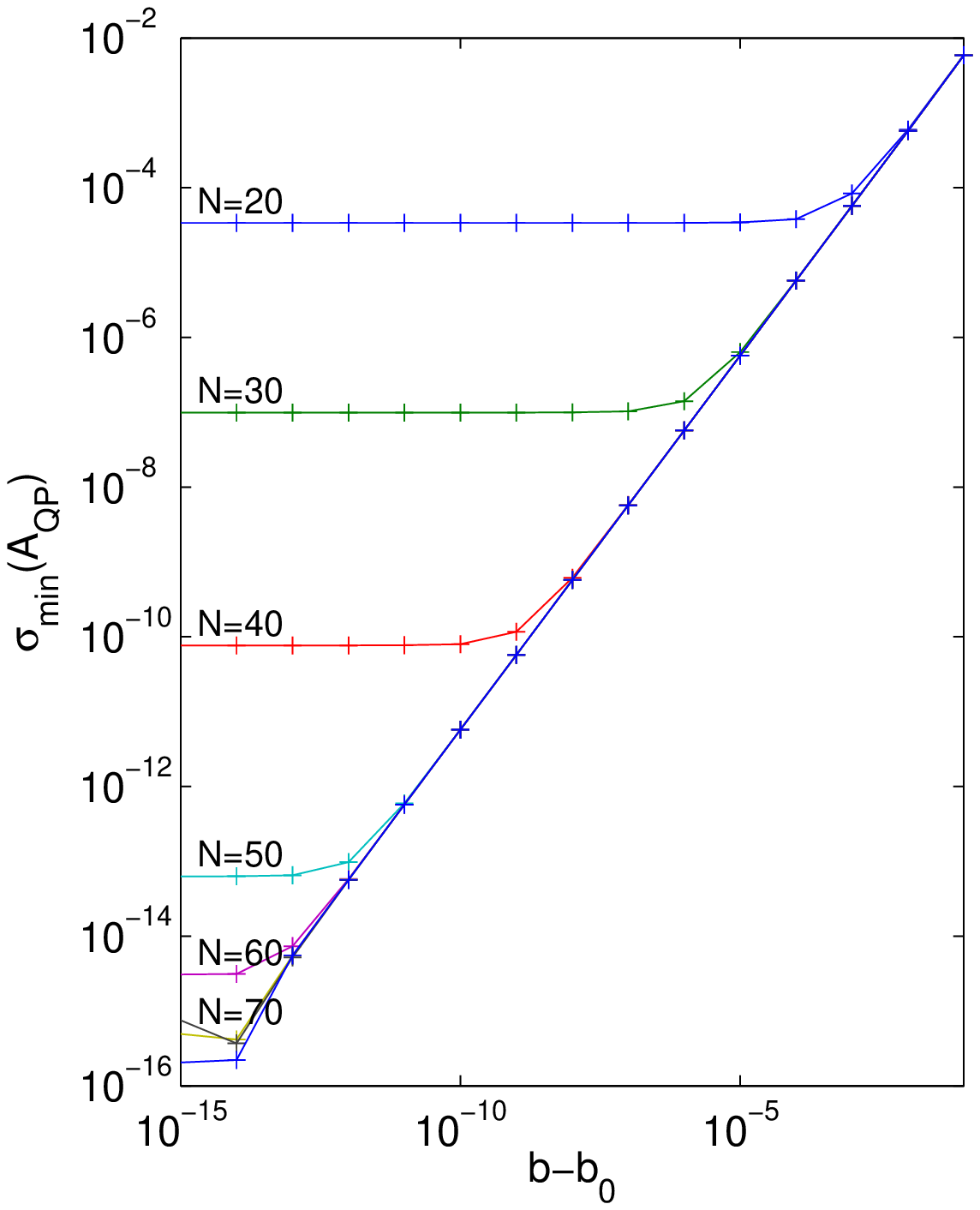}}
\;
c)\raisebox{-2.4in}{\hspace{-2ex}\ig{height=2.3in}{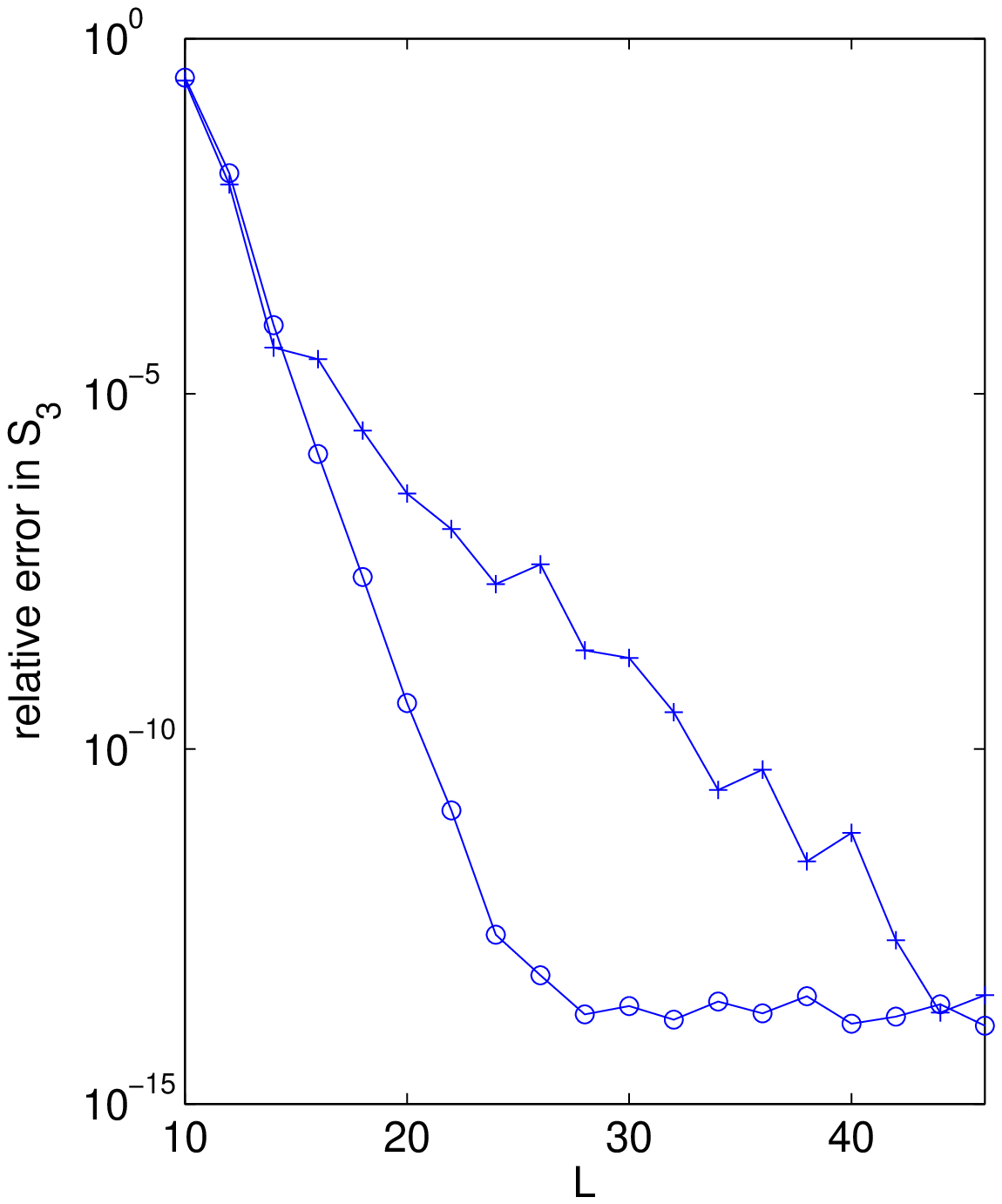}}
}%
\ca{Convergence for quasi-periodic Greens function scheme of Sec.~\ref{s:gqp}.
a) Absolute error in
$\sigma_\tbox{min}(\tilde{A}_\tbox{QP})$
vs $N$ the number of quadrature nodes on $\pO$,
for Bloch parameters $a=\pi/2$ and the two different $b$ values labeled.
The unit cell
with $\eb=(1,0)$, $\el=(0.4,1)$,
and inclusion, described by the radial function
$r(\theta) = 0.2(1 + 0.3 \cos3\theta)$, are shown in the inset.
Index is $n=3$ and frequency $\om=4.5$.
For $b=2$, error is taken relative to the converged value $0.01879908530381247$;
for $b=b_0$, relative to 0.
The matrix $\tilde{A}_\tbox{QP}$ has 2-norm of about 25.
b) $\sigma_\tbox{min}(\tilde{A}_\tbox{QP})$ vs difference in parameter $b$
from the Bloch eigenvalue $b_0$, for several different numbers of
quadrature points $N$. Note the horizontal log scale.
This shows that it is the convergence rate at the Bloch eigenvalue
that controls the accuracy with which the minimum can be found.
c) Relative error ($+$ symbols)
in evaluation of lattice sum $S_3$ by the method of Sec.~\ref{s:eval}
vs the maximum order $L$ in \eqref{e:jexp}.
Parameters are as in Table~II of \cite{McPh00},
whose claim $S_3 = 2.13097899279352 + 5.66537068305984i$
is taken as the true value.
Also, relative error ($\circ$ symbols) for  $\tilde{S}_3$
which excludes 
the $3\times3$ block of neighbors
(parameters are the same; true value is taken as the
converged value at $L=50$).
}{f:gqpconv}
\efi

\subsection{New method for evaluation of the quasi-periodic Greens function}
\label{s:eval}

In order to compute the elements of $\tilde{A}_\tbox{QP}$,
one must evaluate $\gqpr$ defined by \eqref{e:gqpr}; in this section, we present a
surprisingly simple (and apparently new) method for this.
Since the sums \eqref{e:gqp} and \eqref{e:gqpspec} converge
too slowly to be numerically useful, many sophisticated schemes have been devised.
Some of these are based on the Fourier representation (such as 
 \cite{Beylkin08}), but most are based on the observation that
\be
\gqpr(r,\theta) =
\sum_{l=-L}^{L} S_l J_l(\omega r) e^{il\theta} ~,
\label{e:jexp}
\ee
where $(r,\theta)$ are the usual polar coordinates, and $J_l$ the
regular Bessel function of order $l$.
As $L\to\infty$, this expression is uniformly convergent in the unit cell $U$,
as long as
there exists a circle about the origin which contains $\overline{U}$ but
encloses no points in $\Lambda\setminus\{\mbf{0}\}$. 
The coefficients $S_l$ in this expansion are
know as  {\em lattice sums}, given by 
\[
S_l = 
 \sum_{\substack{m,n\in\mathbb{Z} \\
(m,n) \neq (0,0)}} \al^m \bt^n \, H_l^{(1)}(\omega r_{mn}) e^{-il\theta_{mn}} ,
\]
where $(r_{mn},\theta_{mn})$ are the polar coordinates of $m\eb+n\el$,
and $H_l^{(1)}$ is the outgoing Hankel function of order $l$.
Thus, the issue of 
evaluating $\gqpr$ has been reduced to that of tabulating the lattice sums.
This problem itself has a substantial literature (see, for example,
\cite{Chin94,dienst01,leung93,McPh00,Moroz01}). Nevertheless, 
very few papers discuss the problem of empty resonances, at which
point the lattice sums $S_l$ blow up. 
One notable exception is the work of Linton and Thompson \cite{linton07}, who analyze this
blowup for periodic one-dimensional arrays in two
dimensional scattering. They also propose a regularization method to overcome it.

We present here the construction of a small linear system whose solution yields the lattice sums rather easily (away from empty resonances). 
In physical terms, we compute the field induced by the free-space Green's
function $G$, determine how it fails to satisfy quasi-periodicity, and use the representation
\eqref{e:jexp} to enforce quasi-periodicity numerically. More precisely, given a field $u$, 
we define the {\em discrepancy} by
\be
d        
\;=\; \vf{f}{f'}{g}{g'} \;:=\; \left[\begin{array}{l}
u|_L - \alpha^{-1}u|_{L+\eb}\\
\un|_L - \alpha^{-1}\un|_{L+\eb}\\
u|_B - \beta^{-1}u|_{B+\el}\\
\un|_B - \beta^{-1}\un|_{B+\el}
\end{array}\right]~.
\label{e:discrep}
\ee
We can interpret $f$, $f'$ as functions on
wall $L$ and $g$, $g'$ as functions on wall $B$. We
construct a $4M$-component column vector $\mbf{d}$ by
sampling these
four functions at Gaussian quadrature points $\{\yy^{(L)}_m\}_{m=1}^M$ on $L$,
and $\{\yy^{(B)}_m\}_{m=1}^M$ on $B$. If we let
the field $u(\xx) = G(\xx)$, then for $m=1,\ldots,M$, the $m$th element of $\mbf{d}$ 
is $G(\yy^{(L)}_m) - \alpha^{-1} G(\yy^{(L)}_m + \eb)$. The 
remaining $3M$ entries in $\mbf{d}$ are computed in the analogous fashion.

Now let $H$ be a (complex) matrix of size
$4M \times (2L+1)$, defined as follows.
For $l=-L,\ldots,L$, fill the $(l+L+1)$th column 
in the same manner as $\mbf{d}$, but using the field $u(\xx) = 
J_l(\omega r) e^{il\theta}$. Letting $\mbf{s}:=\{S_l\}_{l=-L}^L$, it is straightforward to verify 
that the linear system
\be
H \mbf{s} = -\mbf{d}
\label{e:slsys}
\ee
yields values for the lattice sums that annihilate the discrepancy induced by the source
$G$. We solve the linear system in the least squares sense. 
This has to be done with some care, since the Bessel functions $J_l$
become exponentially small for large $l$. A simple fix is to 
right-precondition the system by scaling the $(l+L+1)$th column of $H$ by
the factor $\rho_l := 1/J_l(\min[\om R,l])$, where $R:=\max_{\xx\in U}|\xx|$
is the unit cell radius. 
The entire procedure may be interpreted as finding the representation
\eqref{e:jexp} which minimizes the $L^2$-norm of
the discrepancy of the resulting $\gqp$.

Fig.~\ref{f:gqpconv}b shows that the error
in evaluating $S_l$, for $l=3$, has exponential convergence in $L$.
We fixed $M=24$ (large enough that further increase had no effect).
14 digits of relative
accuracy are achieved for $L\ge 46$, comparable in accuracy to \cite{McPh00}.
Although the maximum achievable accuracy for $S_l$ deteriorates exponentially
as $|l|$ increases,
the resulting accuracy of $\gqpr$ computed via \eqref{e:jexp} is close to
14 digits everywhere in $U$.%
\footnote{This is to be expected from arguments similar to \cite[Eq.~(5)]{mfs}:
  the residual of the linear system, around $10^{-14}$, approximates
  the boundary error norm, which in turn controls the interior error norm
  when using a basis of particular solutions to the Helmholtz equation.}
We do not claim that our method is optimal in terms of speed
(although at 0.05 sec to solve for all $S_l$ values, it is adequate),
merely that it is accurate, convenient and robust. To our knowledge
it has not been proposed in the literature. 

The convergence rate in the boundary $L^2$-norm
of expansions such as \eqref{e:jexp} depends on the
(conformal) distance from the domain to
the nearest field singularity (a result of Vekua's theory and approximation
in the complex plane \cite[Ch.~6]{timothesis}).
Thus, the rate may be improved by increasing this distance by
removing the rest of the
$3\times 3$ block of nearest neighbors from the
lattice sum, and representing
\be
\gqprr(\xx)\;:= \!\!
\sum_{(j,k)\;\in\;\mathbb{Z}^2 \setminus \{-1,0,1\}^2}
\!\!\!\!\al^j \bt^k G(\xx-j\eb-k\el)
 = \sum_{l=-L}^{L} \tilde{S}_l J_l(\omega r) e^{il\theta} ~.
\label{e:gqp3}
\ee
To solve for $\{\tilde{S}_l\}$, the right-hand side of the linear system is
now chosen to be
the direct summation of these neighbors, $u(\xx) = \tilde{G}(\xx) :=
\sum_{j,k\in \{-1,0,1\}} \al^j \bt^k G(\xx-j\eb-k\el)$.
We may then evaluate $\gqp = \tilde{G} + \gqprr$.
As Fig.~\ref{f:gqpconv}c shows, the convergence rate for $\tilde{S}_l$,
and hence for $\gqp$, is now a factor 2--3 better.
Hence we use this method below, fixing $L=30$.

\bfi
\hspace{-5ex}
\mbox{%
\ig{height=2.1in}{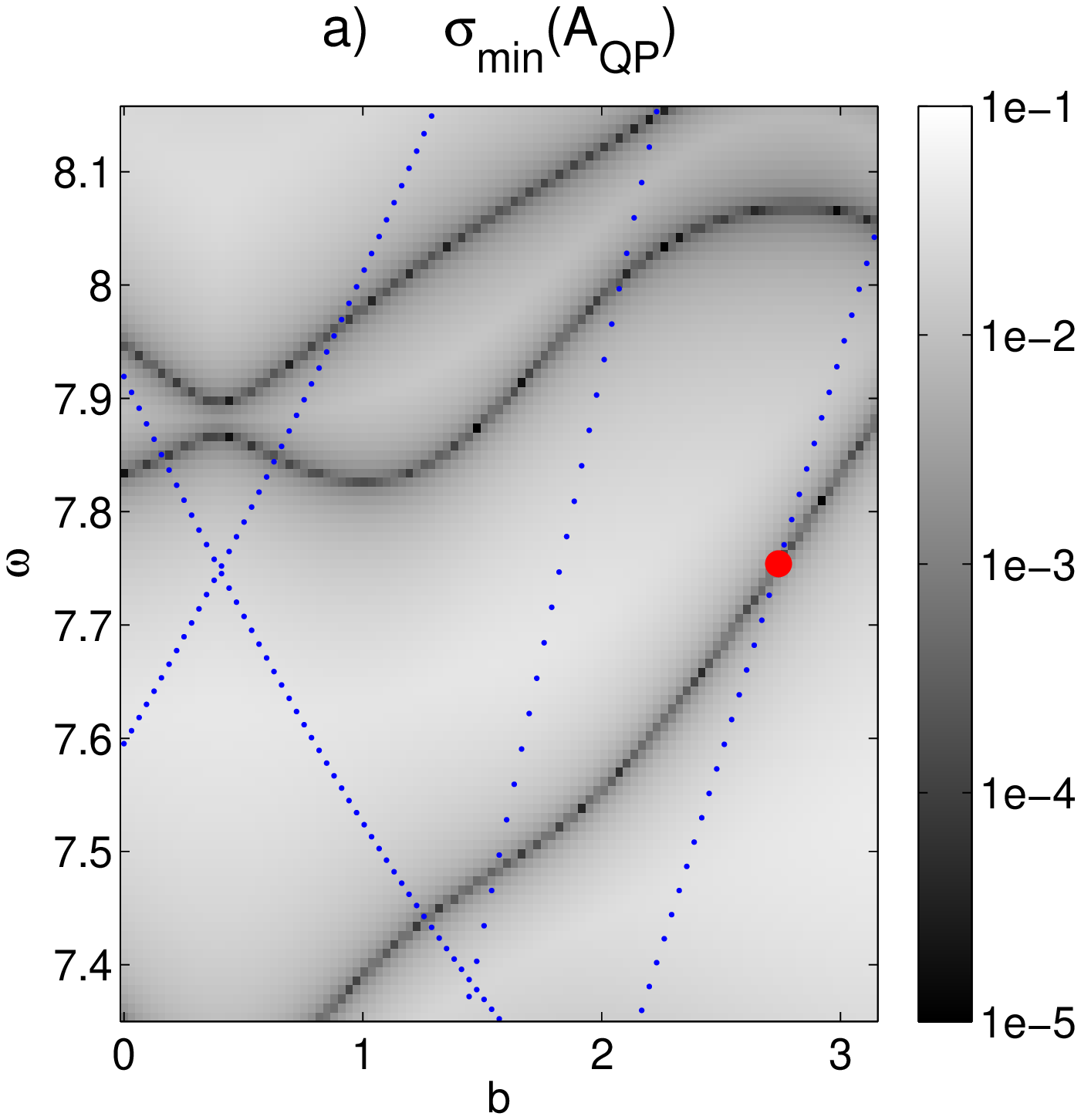}
\ig{height=2.1in}{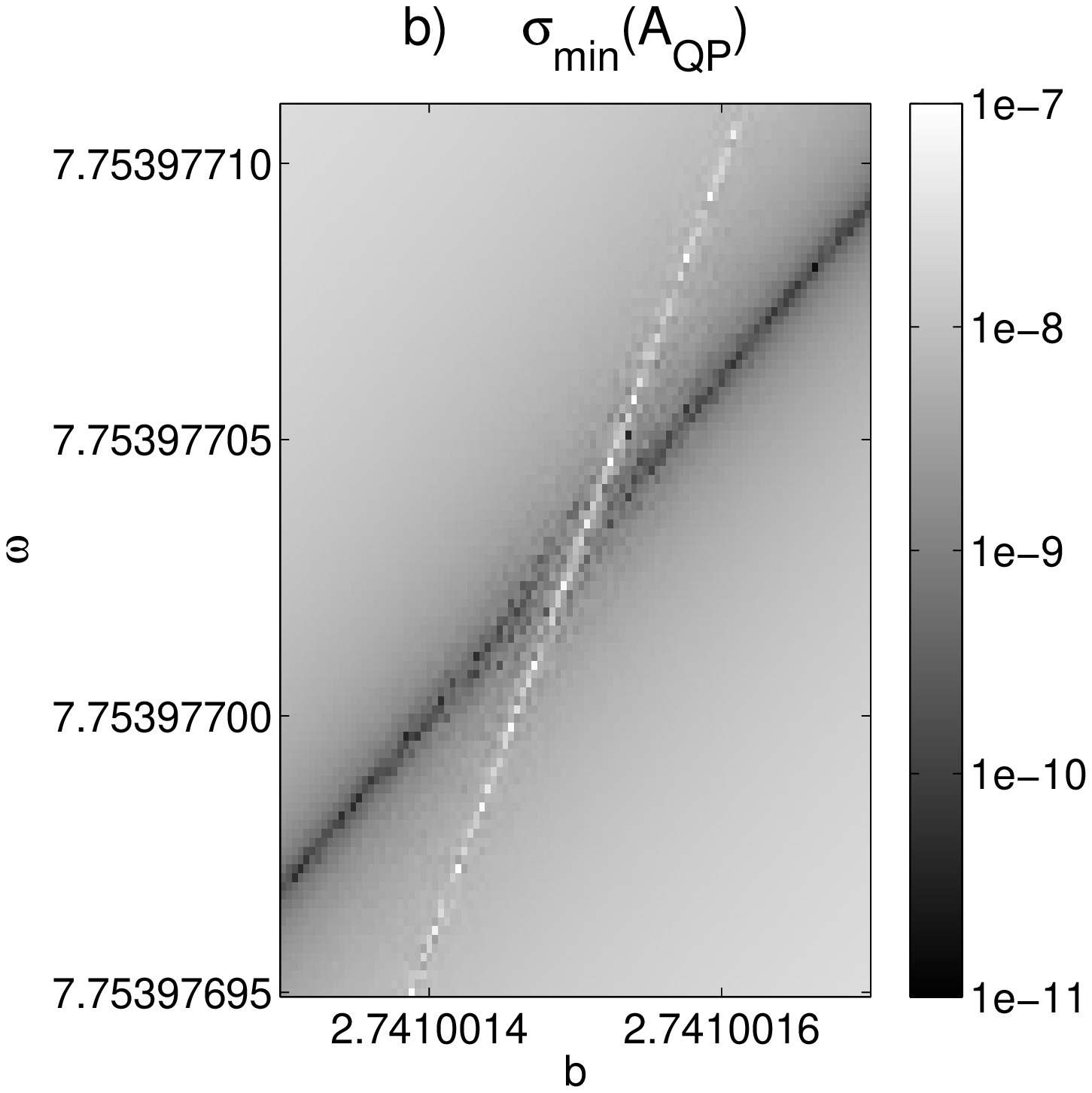}
\ig{height=2.1in}{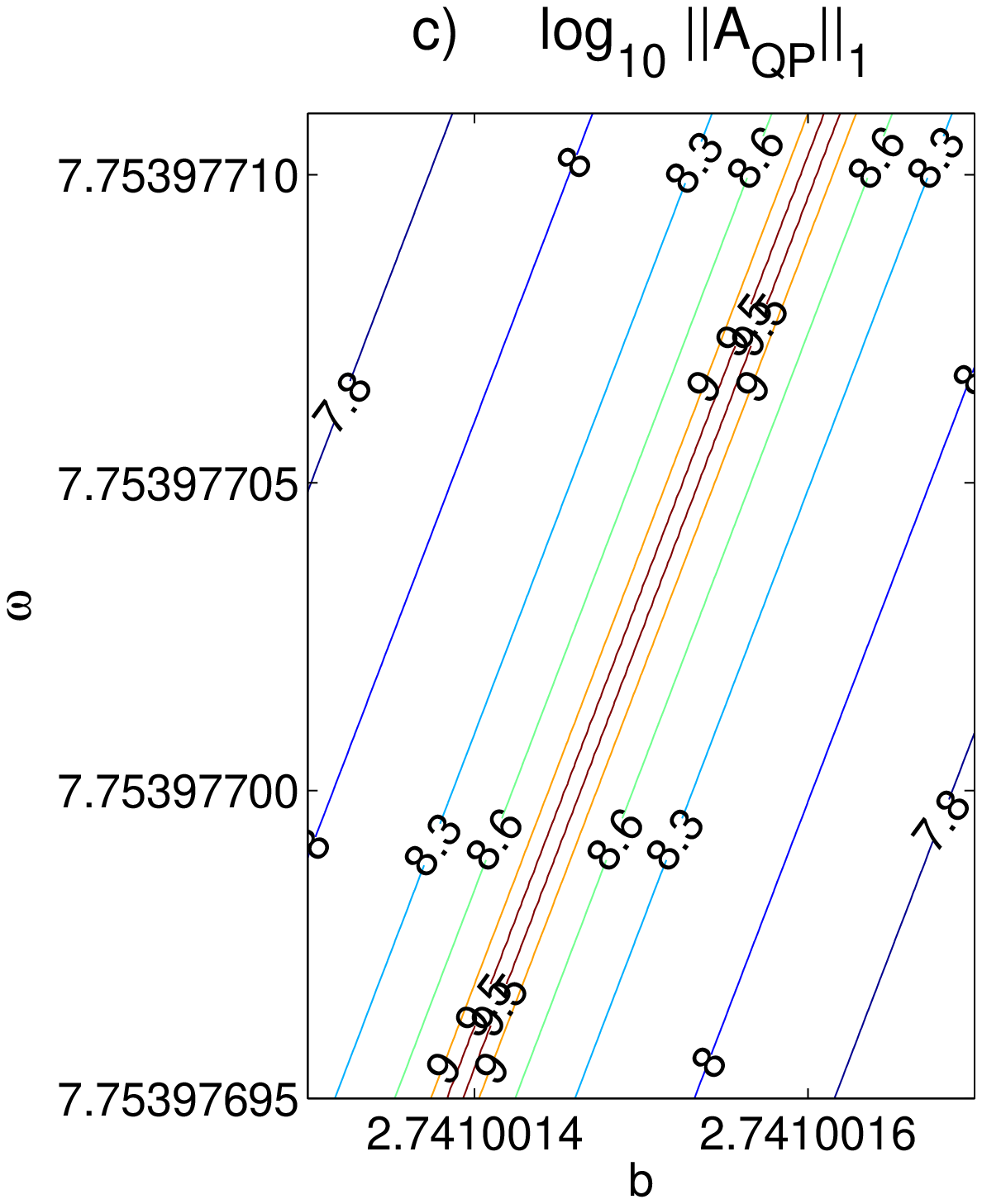}
}%
\ca{Breakdown of quasi-periodic Greens function scheme,
for the system of Fig.~\ref{f:gqpconv}a except with $\el=(0.5,1)$. 
a) Minimum singular value of $\aqpm$ vs $b = \kk \cdot \el$ and $\omega$,
as a log density plot over a slice with fixed $a=0.8$.
Dark curves indicate the band structure, and superimposed
dotted lines the `empty' band structure where $\gqp$ blows up.
b) Zooming in by a factor of $10^7$ to the region shown by the dot in a), 
showing failure to resolve band structure at the intersection.
c) log 1-norm of the matrix $\tilde{A}_\tbox{QP}$ plotted over the same
region as b); it is of order the inverse of the distance to
the empty band structure.
}{f:break}
\efi

\subsection{The empty resonance problem} 
\label{s:flaw}

Given a photonic crystal (inclusion $\Omega$ with index $n$),
using the methods of Sections \ref{s:discr} and \ref{s:eval} we are able to
construct the matrix $\aqpm$ for any given frequency and Bloch parameters
$(\om,a,b)$.
Fig.~\ref{f:break}a shows the minumum singular value of this matrix
as a function over the $(b,\om)$ plane, for constant $a$: the band structure
is visible as the zeros of this function.
We have also superimposed the band structure of the empty unit cell (dotted lines).
Theorem~\ref{t:aqp} guarantees that,
away from the empty unit cell band structure,
no spurious modes will be found, and that no modes are missed.

However, zooming in to one of the many intersections of the two sets of curves
(Fig.~\ref{f:break}b), we see that in the neighborhood of the
empty band structure,
the desired singular values take on arbitrary fluctuating values
that obscure the theoretical behavior near their intersection.
This prevents any attempt to locate the desired zero set to an accuracy
better than $O(\sqrt{\emach})$, where $\emach$ is the machine precision.
As Fig.~\ref{f:break}c shows, this is explained by
the blowup of the entries of the matrix  $\aqpm$ as one approaches
the empty band structure. This, in turn,
causes unbounded roundoff error when computing small singular values
in finite-precision arithmetic.
\begin{rmk}
The above demonstrates a fundamental flaw inherent in the use of
the quasi-periodic Greens function  in
band structure problems; there are empty-resonant parameter sets
(sheets in the space $(\om,a,b)$)
where the desired band structure cannot be computed.
Furthermore, loss of accuracy is inevitable near these parameter sets.
\end{rmk}
This motivates the development of a more robust scheme.

\section{Periodizing using auxiliary densities on the unit cell walls}
\label{s:qplp}


\subsection{Inclusion images and a new linear system}

Section~\ref{s:eval} illustrated the fact, well known in
the fast multipole literature \cite{fmmlattice,CMS,fmm1,fmm2,dienst01},
that summing the nearest neighbors directly
(i.e.\ excluding them from the quasi-periodic field representation)
results in much improved convergence rates.
This motivates defining generalizations of \eqref{e:s} and \eqref{e:d}
that include summation over the appropriately phased $3\times 3$
nearest neighbor images, as shown in Fig.~\ref{f:g}b,
\bea
(\tilde{\cal S}^{(\om)}\sigma)(\xx) &=&
\int_{\pO} \sum_{j,k\in\{-1,0,1\}} \!\!\al^j \bt^k
G^{(\om)}(\xx,\yy+j\eb+k\el) \,\sigma(\yy) ds_\yy
\label{e:st}
\\
(\tilde{\cal D}^{(\om)}\tau)(\xx) &=&
\int_{\pO}
\sum_{j,k\in\{-1,0,1\}} \!\!\al^j \bt^k
\frac{\partial G^{(\om)}}{\partial n_\yy}(\xx,\yy+j\eb+k\el) \,\tau(\yy) ds_\yy
\label{e:dt}
\eea

We now choose a 
layer potential representation for $u$ that involves only free space kernels:
\be
u \;=\; \left\{\begin{array}{ll}
{\cal S}^{(n\om)}\sigma + {\cal D}^{(n\om)}\tau
& \mbox{in } \Omega\\
\tilde{\cal S}^{(\om)}\sigma + \tilde{\cal D}^{(\om)}\tau  + u_\tbox{QP}[\xi]
& \mbox{in }
U\setminus\overline{\Omega}
\end{array}\right.
\label{e:urep}
\ee
The auxiliary field
$u_\tbox{QP}$ will be represented by a new set of layer potentials
that lie on the ``tic-tac-toe'' stencil of Fig.~\ref{f:g}b, consisting of the 
boundary of $U$ and its closest extensions, none lying in the interior of $U$.
We will return to this in section \ref{s:xi}.
For the moment, let us denote the unknown densities that determine
$u_\tbox{QP}$ by $\xi$. By construction, the representation
\eqref{e:urep} satisfies \eqref{e:pdei} and \eqref{e:pdee} in $U$,
so that it remains only to impose both the matching/continuity conditions \eqref{e:cont}
and quasi-periodicity \eqref{e:f}-\eqref{e:gp}. 
Imposing the mismatch $m$ defined by \eqref{e:m} and the 
discrepancy $d$ defined by \eqref{e:discrep} on $u$, the unknowns in 
\eqref{e:urep} must satisfy a linear system of the form:
\be
E \vt{\eta}{\xi} := \mt{A}{B}{C}{Q} \vt{\eta}{\xi} = \vt{m}{d}~,
\label{e:2x2}
\ee
where, as before, $\eta:=[\tau;-\sigma]$, 
We will describe the operators $A$, $B$, $C$, and $Q$ in more detail shortly.
For the moment, note that if there exists a density $[\eta;\xi]$
which generates a nontrivial field with vanishing mismatch and discrepancy,
then it is 
a solution to \eqref{e:pdei}-\eqref{e:cont} and \eqref{e:f}-\eqref{e:gp} and
the corresponding parameters $(\om,a,b)$ must be a Bloch eigenvalue.
Numerical evidence supports the following stronger claim, the analog of
Theorem~\ref{t:aqp}.
\begin{cnj}
$(\om,a,b)$ is a Bloch eigenvalue if and only if $\Null E \neq \{0\}$.
\label{c:qplp}
\end{cnj}
This suggests, as in Section~\ref{s:gqp},
computing the band structure by locating the parameter families
where (a discretization of) $E$ is singular.
The point of the new scheme is that it should be robust;
if the conjecture
holds, then (in contrast to the quasiperiodic Green's function approach), 
there will be no spurious parameter values where the method breaks down.

To discuss the operators in $E$, we need some additional notation. We assume
that the wavenumber $\om$ and quasiperiodicity parameters $(a,b)$ are given.
Let $W$ be a curve in $\RR$ on which
single and double layer densities are defined,
with the corresponding potentials written as
\bea
({\cal S}_{W} \sigma)(\xx) &=& \int_{W}
G(\xx,\yy) \sigma(\yy) ds_\yy 
\label{e:sw}
\\
({\cal D}_{W} \tau)(\xx) &=& \int_{W}
\frac{\partial G}{\partial n_\yy} (\xx,\yy) \tau(\yy) ds_\yy  \, .
\label{e:dw}
\eea
Letting $V$ be a (possibly distinct) target
curve in $\RR$, we define the operators
\bea
({S}_{V,W} \sigma)(\xx) &=& \int_{W}
G(\xx,\yy) \sigma(\yy) ds_\yy \qquad \xx \in V
\label{e:svw}
\\
({D}_{V,W} \tau)(\xx) &=& \int_{W}
\frac{\partial G}{\partial n_\yy} (\xx,\yy) \tau(\yy) ds_\yy  \qquad\xx \in V
\label{e:dvw}
\\
({D}^\ast_{V,W} \sigma)(\xx) &=& \int_{W}
\frac{\partial G}{\partial n_\xx} (\xx,\yy) \sigma(\yy) ds_\yy  \qquad \xx \in V
\label{e:dvwast}
\\
({T}_{V,W} \tau)(\xx) &=& \int_{W}
\frac{\partial^2 G}{\partial n_\xx\partial n_\yy}(\xx,\yy)
\tau(\yy) ds_\yy  \qquad \xx \in V~.
\label{e:tvw}
\eea
When $V=W$, these operators are to be understood in the principal value sense.
By analogy with 
\eqref{e:st}, \eqref{e:dt}, versions of these operators whose kernels include the 
phased sum over $3\times 3$ images of the source
are indicated with a tilde ($\sim$): that is,
$\tilde{S}_{V,W}, \tilde{D}_{V,W},\tilde{D}^\ast_{V,W}$, and $\tilde{T}_{V,W}$.

We are now in a position to provide explicit expressions for the operators 
$A, B, C, Q$ in \eqref{e:2x2}. Comparing \eqref{e:urep} to \eqref{e:qprep}, it is clear that
the operator $A$ is the same as $\aqp$ in \eqref{e:aqp} but with the
replacement of
$S_\tbox{QP}^{(\om)}$,  $D_\tbox{QP}^{(\om)}$ and $T_\tbox{QP}^{(\om)}$,
by $\tilde{S}_{\pO,\pO}$, $\tilde{D}_{\pO,\pO}$ and $\tilde{T}_{\pO,\pO}$,
respectively.
It is straightforward to verify that $A$ is a compact perturbation of the identity.

The operator $C$ describes the effect of the inclusion densities on
the discrepancy $d$.
Its eight sub-blocks are found by inserting \eqref{e:st} and \eqref{e:dt}
into \eqref{e:urep} then evaluating \eqref{e:discrep}, giving
\[
C \; =\;
\begin{bmatrix}
\tilde{D}_{L,\pO}\!-\al^{-1}\tilde{D}_{L+\eb,\pO}&
-\tilde{S}_{L,\pO}\!+\al^{-1}\tilde{S}_{L+\eb,\pO}\\
\tilde{T}_{L,\pO}\!-\al^{-1}\tilde{T}_{L+\eb,\pO}&
-\tilde{D}^\ast_{L,\pO}\!+\al^{-1}\tilde{D}^\ast_{L+\eb,\pO}\\
\tilde{D}_{B,\pO}\!-\bt^{-1}\tilde{D}_{B+\el,\pO}&
-\tilde{S}_{B,\pO}\!+\bt^{-1}\tilde{S}_{B+\el,\pO}\\
\tilde{T}_{B,\pO}\!-\bt^{-1}\tilde{T}_{B+\el,\pO}&
-\tilde{D}^\ast_{B,\pO}\!+\bt^{-1}\tilde{D}^\ast_{B+\el,\pO}
\end{bmatrix}
\]
Consider now the any of the four upper sub-blocks of $C$. There are nine phased copies of 
$\pO$ which contribute to the field on the left ($L$) and right
($L+\eb$) wall.
From symmetry and translation invariance considerations, however,
it is easy to check that the contributions from the six left-most images on $L$ (dotted curves
in Fig.~\ref{f:cancel}a) are equal to the contributions of the six right-most images
on $L+\eb$ (dotted curves in Fig.~\ref{f:cancel}b). In the $(1,1)$ sub-block, for example,
we have:
\[
\tilde{D}_{L,\pO}-\al^{-1}\tilde{D}_{L+\eb,\pO}
\; =\;
\sum_{k\in\{-1,0,1\}} \bt^k \left(
\al D_{L,\pO+\eb+k\el} - \al^{-2} D_{L,\pO-2\eb+k\el}
\right)
\]
A rotated version of the analysis applies to the 
lower four sub-blocks in $C$.
The result is that the entries in $C$ involve only source-target interactions at
distances {\em greater} than the size of the unit cell,
ensuring the rapid convergence of a representation in terms of smooth functions.

We next discuss the representation of $\xi$ and $u_\tbox{\rm QP}[\xi]$ in more detail,
which will determine the form of blocks $Q$ and $B$ of the full system matrix $E$.

\bfi 
\centering\ig{width=0.9\linewidth}{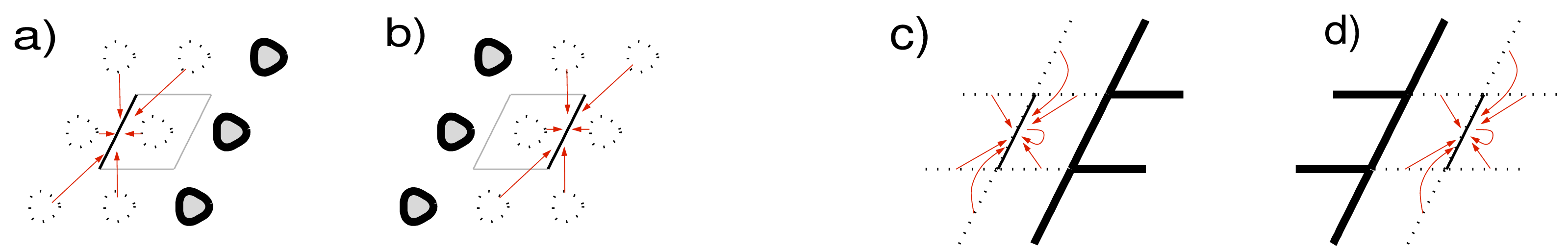}
\ca{Discrepancy cancellation
due to neighbor image sums.
Each arrow represents the influence of a source density on a target segment.
a) For the four upper sub-blocks of $C$, the six nearest source images (dotted) contribute
to the discrepancy on the left wall $L$. b) The six nearest source images (dotted)
contribute exactly the same field (suitable phased) to the right wall $(L + \eb)$. 
The net result is that only the distant sources, shown in bold, have a non-zero effect.
The same holds for all four upper sub-blocks of $C$. A rotated version applies
to the lower four sub-blocks of $C$.
c,d) Contributions to the  sub-blocks $Q_{LL}$ and $Q_{LB}$ of $Q$. The 
seven indicated terms (dotted source segments) cancel in the two diagrams,
leaving only the contributions from distant wall segments shown in bold.
A rotated version applies to the sub-blocks $Q_{BL}$ and $Q_{BB}$.
}{f:cancel}
\efi 

\subsection{Choice of auxiliary densities and their images}
\label{s:xi}

The auxiliary field $u_\tbox{QP}$ is determined by the choice of layer potentials
on the boundary of (and outside of) $U$. We will use double and single layer 
densities
on both the left $(L)$ and bottom $(B)$ boundaries of $U$,
as well as on the other segments of the ``tic-tac-toe'' board in
Fig.~\ref{f:g}b. More precisely, we 
define the vector of unknowns $\xi$ by $\xi := [\tau_L;-\sigma_L;\tau_B;-\sigma_B]$,
and set
\be
u_\tbox{QP}  = \!\!\! \sum_{\substack{j\in\{0,1\}\\ k\in\{-1,0,1\}}}
\!\!\! \al^j \bt^k
\left( {\cal S}_{L+j\eb+k\el}\sigma_L +
  {\cal D}_{L+j\eb+k\el}\tau_L \right)
\;\;+
\!\!
\sum_{\substack{j\in\{-1,0,1\} \\ k\in\{0,1\}}}
\!\!\! \al^j \bt^k
\left( {\cal S}_{B+j\eb+k\el}\sigma_B +
  {\cal D}_{B+j\eb+k\el}\tau_B \right)
\label{e:uqp}
\ee
The inclusion of the image segments leads to cancellations that are numerically
advantageous in the operator $Q$, just as we found that images helped with
the operator $C$ in the preceding section.

We should first clarify the definition \eqref{e:discrep} of the discrepancy
functions: field values should be interpreted as their limiting values
on the wall approaching from {\em inside} $U$, since it is
the field in $U$ that \eqref{e:urep} and \eqref{e:uqp} represent.
For example, $f : = u^+|_L - \al^{-1} u^-|_{L+\eb}$, where, as before
$u^{\pm}(\xx) : = \lim_{\eps\to0^+}
u(\xx\pm\eps\mbf{n})$, and $\mbf{n}$ is the normal at $\xx$.

Recall now that the operator $Q$ expresses the effect of the four densities in $\xi$
on the four discrepancy functions $f$, $f'$, $g$, $g'$.
If \eqref{e:uqp} contained only the terms $j=k=0$, this would
correspond to densities $\sigma_L$ and $\tau_L$ placed on $L$,
and $\sigma_B$ and $\tau_B$ placed on $B$.
While this is mathematically acceptable, it results in various complicated
self-interactions and interactions between segments that share a common corner.
This would lead to singularities in densities requiring more complicated
discretization and quadrature. Although there has been significant progress in 
this direction (see, for example, \cite{bremer,helsing_corner}), in the present 
context we have the luxury of including the ten additional image segments in
\eqref{e:uqp}, which cancel both the self and near-field corner interactions.
As a result,
our implementation is simpler and involves fewer degrees of freedom.
The cancellation mechanism is shown in Fig.~\ref{f:cancel}.
The effect on $u^+|_L$ of the seven segments touching $L$, for example,
cancels the effect on $u^-|_{L+\eb}$ of the seven segments touching $L+\eb$,
leaving only ten far field contributions.

It is important to note that  the local terms due to the jump relations do not cancel:
e.g.\ a density function $\tau_L$ placed on $L$ contributes a term
$\frac{1}{2} \tau_L$ to $u^+|_L$, while
$\al\tau_L$ placed on $L+\eb$ contributes $-\frac{1}{2} \al \tau_L$
to $u^-|_{L+\eb}$. These two terms add to contribute $\tau_L$ to $f$.
One may check in this fashion that the jump relations contribute
an identity to the diagonal sub-blocks of $Q$.
This yields the crucial result that $Q$ is the identity plus a compact operator,
with the compact part generated by interactions at a distance greater than 
the size of the unit cell.
After the above cancellations and simplification, we have,
\[
Q = I \,+\,
\begin{bmatrix}
Q_{LL} & Q_{LB} \\
Q_{BL} & Q_{BB} 
\end{bmatrix}
\]
where
\begin{align*} 
Q_{LL} &= 
\begin{bmatrix}
\displaystyle
\sum\limits_{j\in\{-1,1\},k\in\{-1,0,1\}} \!\!\! j \al^j\bt^k
D_{L,L+j\eb+k\el} &
\qquad - \hspace{-3ex}\displaystyle
\sum\limits_{j\in\{-1,1\}, k\in\{-1,0,1\}} \!\!\!
j \al^j\bt^k S_{L,L+j\eb+k\el}
\\ \\
\displaystyle
\sum\limits_{j\in\{-1,1\},k\in\{-1,0,1\}} \!\!\! j \al^j\bt^k
T_{L,L+j\eb+k\el} &
\qquad - \hspace{-3ex}\displaystyle
\sum\limits_{j\in\{-1,1\},k\in\{-1,0,1\}} \!\!\!
j \al^j\bt^k D^\ast_{L,L+j\eb+k\el}
\end{bmatrix}
\\ \\
Q_{LB} &= 
\begin{bmatrix}
\displaystyle
  \sum\limits_{k\in\{0,1\}} \!\bt^k\bigl(\al D_{L,B+\eb+k\el}-\al^{-2}D_{L,B-2\eb+k\el}\bigr)
&
 \displaystyle
  \sum\limits_{k\in\{0,1\}} \!\bt^k\bigl(-\al S_{L,B+\eb+k\el}
+\al^{-2}S_{L,B-2\eb+k\el}\bigr)
\\ \\                                     
 \displaystyle
  \sum\limits_{k\in\{0,1\}} \!\bt^k\bigl(\al T_{L,B+\eb+k\el}
-\al^{-2}T_{L,B-2\eb+k\el}\bigr)
&
 \displaystyle
  \sum\limits_{k\in\{0,1\}} \!\bt^k\bigl(-\al D^\ast_{L,B+\eb+k\el}
+\al^{-2}D^\ast_{L,B-2\eb+k\el}\bigr)
\end{bmatrix}
\end{align*}
\begin{align*}
Q_{BL} &= 
\begin{bmatrix}
 \displaystyle
  \sum\limits_{j\in\{0,1\}} \!\al^j\bigl(\bt D_{B,L+j\eb+\el}
-\bt^{-2}D_{B,L+j\eb-2\el}\bigr)
&
 \displaystyle
  \sum\limits_{j\in\{0,1\}} \!\al^j\bigl(-\bt S_{B,L+j\eb+\el}
+\bt^{-2}S_{B,L+j\eb-2\el}\bigr)
\\ \\
 \displaystyle
  \sum\limits_{j\in\{0,1\}} \!\al^j\bigl(\bt T_{B,L+j\eb+\el}
-\bt^{-2}T_{B,L+j\eb-2\el}\bigr)
&
 \displaystyle
  \sum\limits_{j\in\{0,1\}} \!\al^j\bigl(-\bt D^\ast_{B,L+j\eb+\el}
+\bt^{-2}D^\ast_{B,L+j\eb-2\el}\bigr)
\end{bmatrix}
\\ \\
Q_{BB} &= 
\begin{bmatrix}
 \displaystyle
\sum\limits_{j\in\{-1,0,1\},k\in\{-1,1\}} \!\!\! k \al^j\bt^k
D_{B,B+j\eb+k\el}  &
 \displaystyle
\qquad - \hspace{-3ex}
\sum\limits_{j\in\{-1,0,1\},k\in\{-1,1\}} \!\!\! k \al^j\bt^k
S_{B,B+j\eb+k\el} 
\\ \\
\displaystyle
\sum\limits_{j\in\{-1,0,1\},k\in\{-1,1\}} \!\!\! k \al^j\bt^k
T_{B,B+j\eb+k\el}  &
 \displaystyle
\qquad - \hspace{-3ex}
\sum\limits_{j\in\{-1,0,1\},k\in\{-1,1\}} \!\!\! k \al^j\bt^k
D^\ast_{B,B+j\eb+k\el} 
\end{bmatrix}
\end{align*}

Finally, we discuss the $B$ operator from \eqref{e:2x2},
which describes the effect of the auxiliary densities $\xi$ on
the mismatch. As with $A$, since the mismatch involves values on only a single
curve $\pO$, there is no opportunity for cancellation.
Inserting \eqref{e:uqp} into \eqref{e:m} we get
\begin{align}
B \, = \!\!
&\sum_{j\in\{0,1\},k\in\{-1,0,1\}} \!\!\!\! \al^j \bt^k
\begin{bmatrix}
D_{\pO,L+j\eb+k\el} & -S_{\pO,L+j\eb+k\el} & 0&0\\
T_{\pO,L+j\eb+k\el} & -D^\ast_{\pO,L+j\eb+k\el} & 0&0\\
\end{bmatrix}
\;\;+\!    \nonumber  \\
&\sum_{j\in\{-1,0,1\}, k\in\{0,1\}} \!\!\!\! \al^j \bt^k
\begin{bmatrix}
0&0&D_{\pO,B+j\eb+k\el} & -S_{\pO,B+j\eb+k\el}\\
0&0&T_{\pO,B+j\eb+k\el} & -D^\ast_{\pO,B+j\eb+k\el}\\
\end{bmatrix}
\label{e:b}
\end{align}

Summarizing the above, $E$ is a compact perturbation of the identity.
Its blocks
$C$ and $Q$ involve interaction distances greater than the unit cell size.
Its block $A$ involves distances controlled by the shape of the inclusion
and its nearest approach to its neighboring images.
Its block $B$ involves distances determined by the nearest approach of
$\pO$ to $\partial U$.

\subsection{Numerical implementation and discretization of $B$}

We discretize the four blocks of
the integral operator $E$ in \eqref{e:2x2}
to give the matrix $\tilde{E} \in \mathbb{C}^{(2N+4M) \times (2N+4M)}$
as follows. 
We sample the densities on $\pO$ at equispaced points with respect to the 
given definition of the curve, as in Section~\ref{s:discr}.
We sample the densities on the 
walls $L$ and $B$ at $M$ standard Gaussian nodes,
as in Section~\ref{s:eval}.
$A$ is then discretized in the same way as $\aqp$ in Section~\ref{s:discr}
with a mix of the periodic trapezoidal rule and 
Kress' singular quadratures for the self-interaction of $\pO$.
The (Nystr\"{o}m) method \eqref{e:nyst} 
may be used for the off-diagonal block $C$,
and also for the wall's self-interaction $Q$.
No special singular quadratures are needed in $Q$,
due to the cancellations discussed above.

The $B$ operator \eqref{e:b} involves computing the field due to source densities on walls
$L$ and $B$ (and their images shown in Fig.~\ref{f:g}b)
at targets on $\pO$. When the distance from the inclusion to boundary
$\dist(\pO,\partial U)$ is large,
the plain Nystr\"{o}m method may be used to construct the discretized matrix
$\hat{B}$.
We will refer to this as discretization method B1. 
With nodes $\yy_m$ and weights $w_m$ on wall $L$,
and nodes $\xx_j$ on $\pO$, for example,
the term $S_{\pO,L}$ in the (1,2)-block of
\eqref{e:b} becomes the matrix $\hat{S}\in\mathbb{C}^{N\times M}$
with elements
$\hat{S}_{jm} = \frac{i}{4} H^{(1)}_0 (\om |\xx_j - \yy_m|) w_m$.

When $\dist(\pO,\partial U)$
becomes small, of course,  the convergence
rate of method B1 will become unacceptably poor.
However, by construction, for a Bloch eigenfunction
the field \eqref{e:uqp} generated by the
wall densities in $\xi$ has no singularities in the $3\times 3$
neighboring block of unit cells.
Hence these densities remain smooth,
poor convergence being merely due to
inaccurate evaluation of their field close to the walls.
This leaves room for a large number of options:
\ben
\item[B2)] For the rows of $\hat{B}$ corresponding
to target points on $\pO$ that are distance $d_0$ or closer to $\pU$,
use adaptive Gauss-Kronrod quadrature%
  \footnote{This was implemented with \matlab's {\tt quadgk},
    which uses a pair of 15th and 7th order formulae,
    with relative tolerance set to $10^{-12}$.}
with integrand given by the product of the kernel function
and the Lagrange polynomial interpolant \cite[Sec.~8.1]{na}
for the density at the $M$ quadrature points.
$d_0$ is some $O(1)$ constant.
For the other rows, use method B1.
\item[B3)] Project onto an order-$L$ cylindrical $J$-expansion at the origin.
This is done by
computing a representation \eqref{e:jexp} for each of the point monopole or
dipole sources in the quadrature approximation to the source densities on
the walls, and then evaluating this at the target quadrature
points on $\pO$ to fill the elements of $\hat{B}$.
The example term discussed for B1 gives $\hat{S} = R P$,
where the ``source-to-local'' matrix $P\in\mathbb{C}^{(2L+1)\times M}$
has elements
\[
P_{lm} = \frac{i}{4} H^{(1)}_l (\om|\yy_m|) e^{-il\theta_m} w_m
\]
and converts single layer density values to $J$-expansion coefficients.
This follows from Graf's addition formula \cite[Eq. 9.1.79]{a+s}.
The expansion matrix $R\in\mathbb{C}^{N\times(2L+1)}$
has elements $R_{jl} = J_l(\om|\xx_j|) e^{il\phi_j}$.
In the above $\theta_m, \phi_j$ are polar angles of points $\yy_m, \xx_j$
respectively.
Similar formulae apply for double layers and evaluation of derivatives.
To reduce dynamic range (hence roundoff error) we in fact scale the
$J$-expansion by the factors $\rho_l$ of Section~\ref{s:eval}
(this does not change the mathematical definition of $\hat{S}$.)
\item[B4)] Use a more sophisticated quadrature approach, such as those of
\cite{beale,helsing_close,mayo}.
\een
Methods B2-B4 evaluate $u_\tbox{QP}$ due to a spectral interpolant of
the discretized wall densities,
with an accuracy that persists up to the boundary of $U$.
Note that this does not increase the number $M$ of degrees of
freedom associated with each such density.
Since the underlying density is smooth (in fact analytic),
the convergence rate is high and
we are able to keep $M$ very modest.






We have implemented methods B1, B2 and B3.
We use the quadrature weights to scale the matrix $\hat{E}$ to give
$\tilde{E}$ in an analogous fashion to
\eqref{e:aqpw}, so that singular values of $\tilde{E}$
approximate those of $E$. 

Finally, there are many possible ways to locate parameter values $(\om,a,b)$
where $\tilde{E}$ is singular.
In this paper, we will simply plot its smallest singular value
$\sigma_\tbox{min}(\tilde{E})$ vs the Bloch parameters,
as in Section \ref{s:gqp}.

\bfi 
\hspace{-5ex}
\mbox{%
a)\raisebox{-2.3in}{\hspace{-2ex}\ig{height=2.2in}{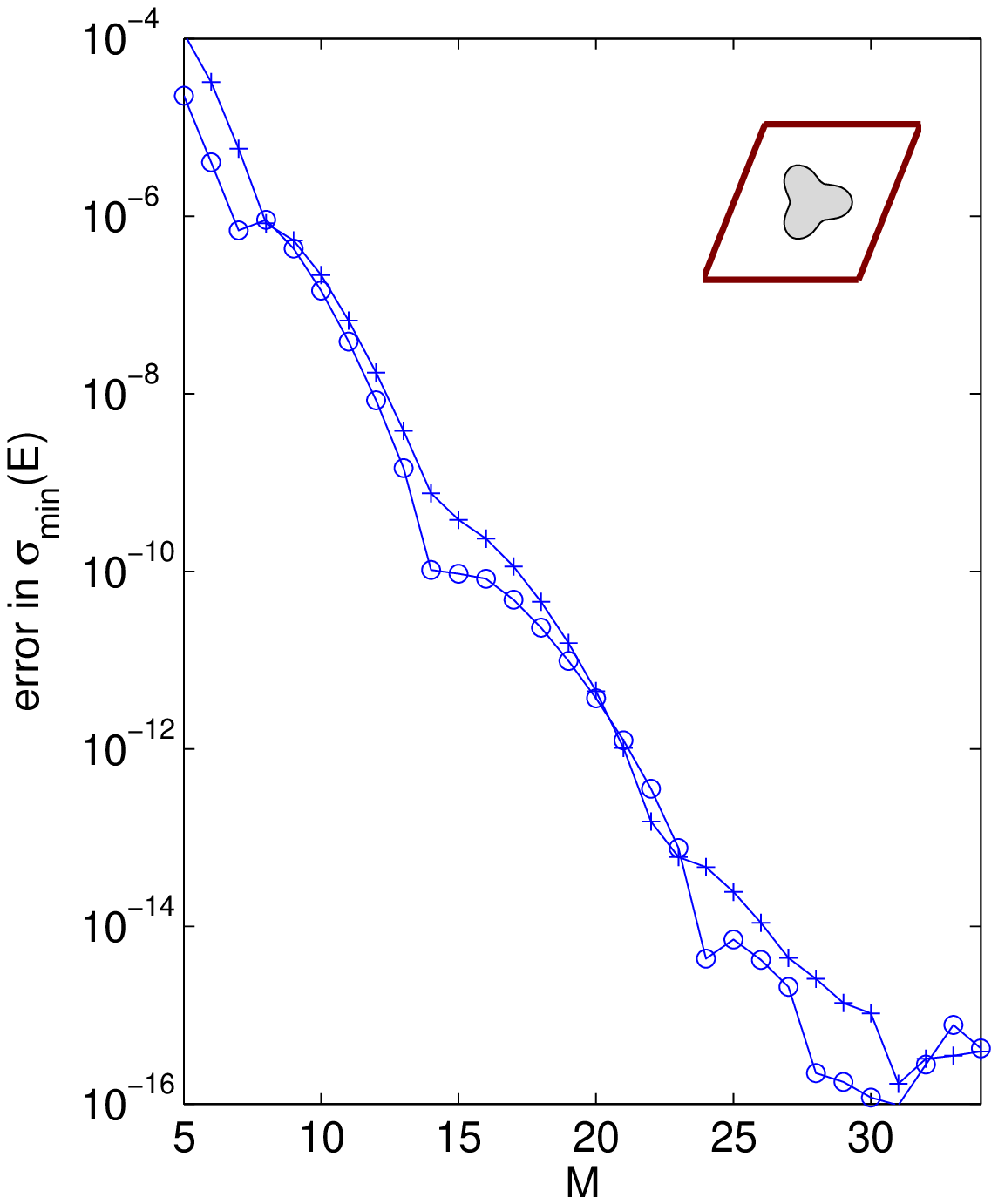}}
b)\raisebox{-2.3in}{\hspace{-2ex}\ig{height=2.2in}{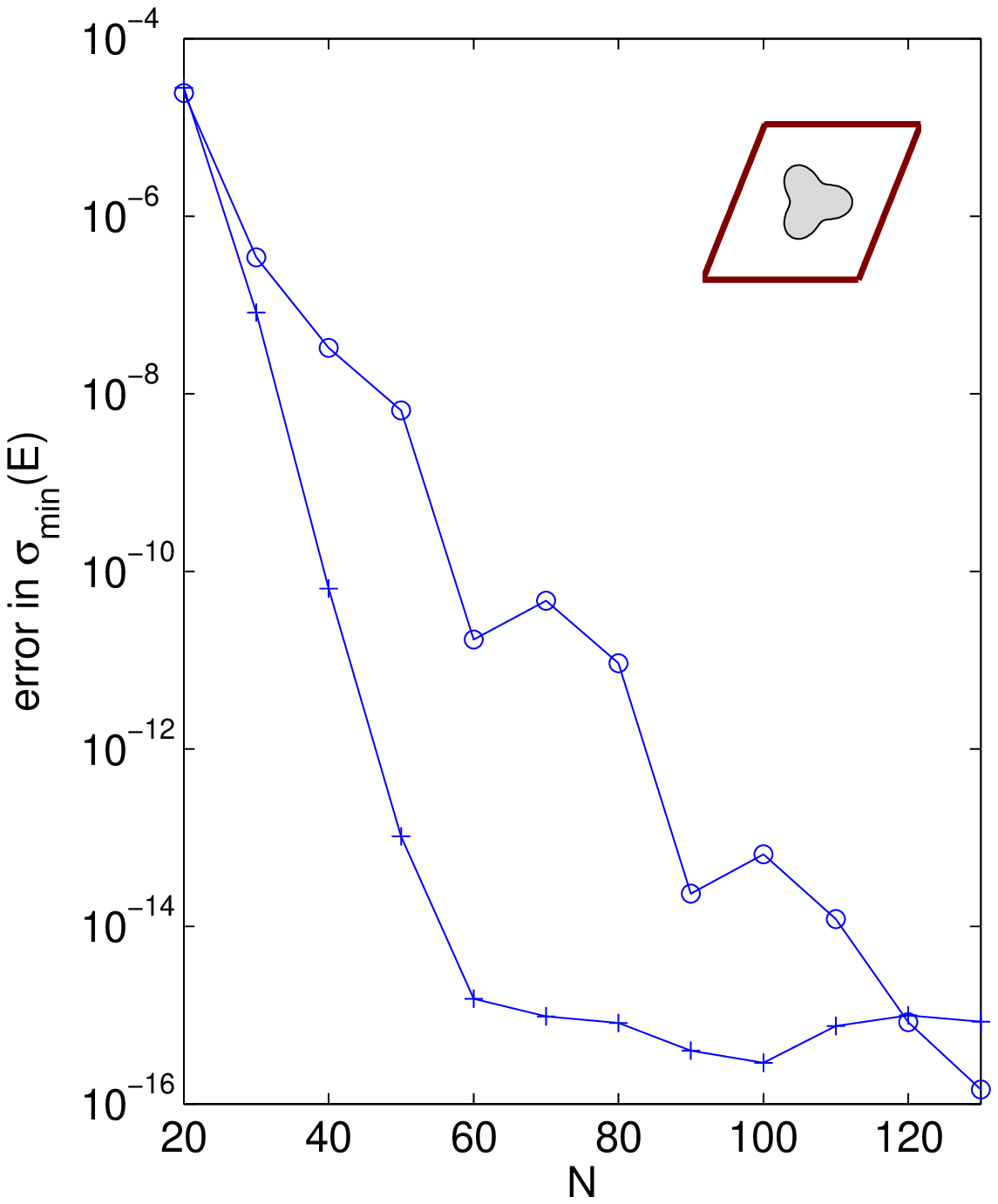}}
c)\raisebox{-2.3in}{\hspace{-2ex}\ig{height=2.1in}{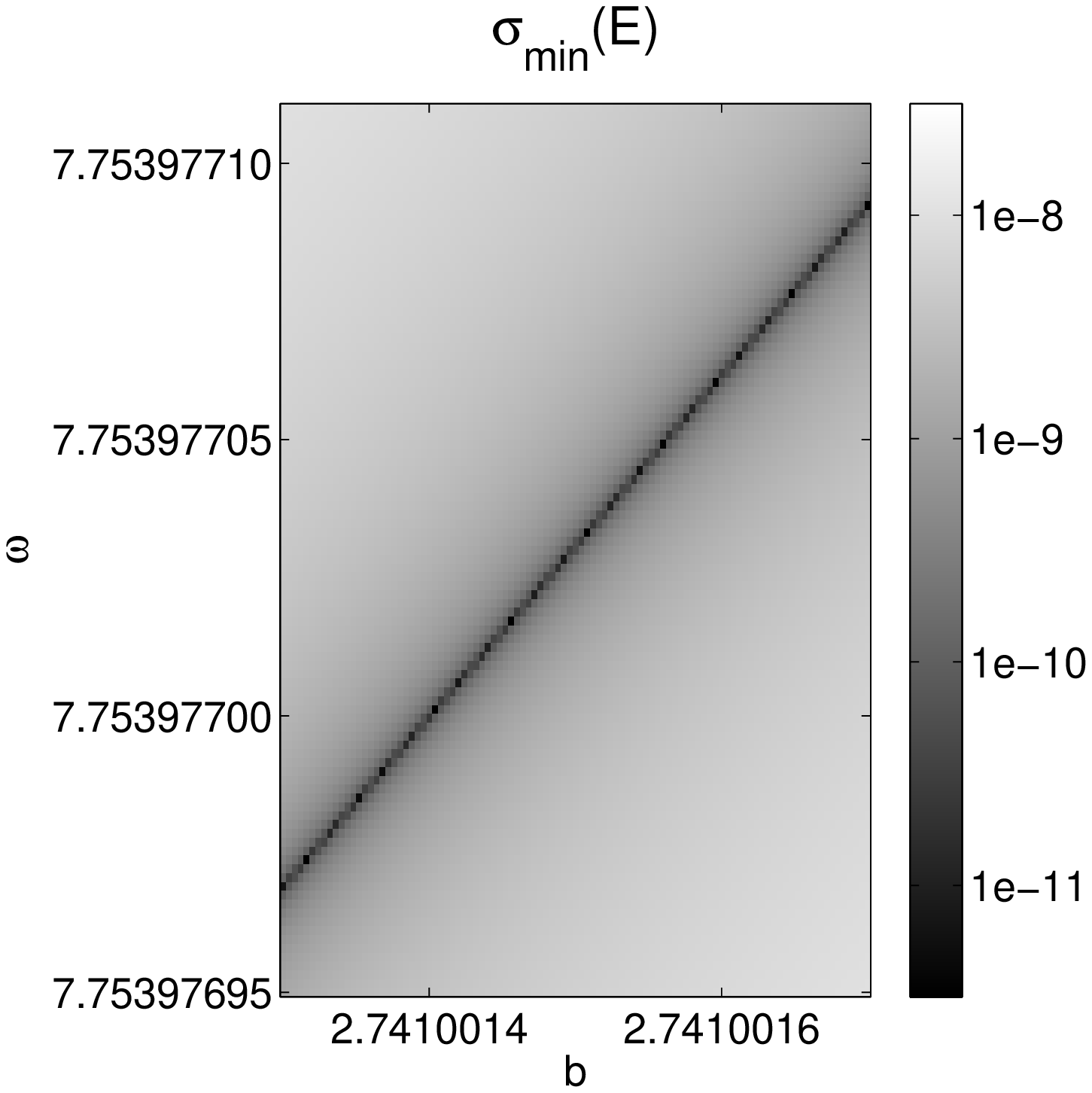}}
}%
\ca{Convergence of new periodizing scheme using auxilliary densities (described
in Section~\ref{s:qplp}), using method B1, for the same geometry and parameters
as in Fig.~\ref{f:gqpconv}a. The meaning of the two curves is also the same as
in the earlier figure. a) Absolute error in $\sigma_\tbox{min}(\tilde{E})$ vs $M$ the
number of nodes on each unit cell wall,
for fixed $N=70$ nodes on $\pO$.
b)
Same as a) except convergence vs $N$, for fixed $M=30$. 
c) Same as Fig.~\ref{f:break}b but using the new scheme:
note the absence of pollution by the empty band structure.
}{f:conv}
\efi

\section{Results of proposed scheme}
\label{s:res}

We first test the convergence of the new scheme for the same
small inclusion used in Section~\ref{s:gqp}, with the simplest discretization 
method for $B$, namely B1.
As before, we test two Bloch parameter $b$ values, one which is
far from an eigenvalue, and one of which is guaranteed to be an eigenvalue
according to Theorem~\ref{t:aqp}.
Fixing $N=70$, which was found in Section~\ref{s:discr} to be
fully converged when at an eigenvalue,
we first vary $M$, the number of nodes per unit cell wall.
Fig.~\ref{f:conv}a shows the convergence of the minimum singular value
of the discretized matrix $\tilde{E}$ to its
converged value (when far from an eigenvalue), or to zero (when at an
eigenvalue). The convergence is spectral, and in both cases full machine accuracy
is reached at $M=30$.
(For $N>70$ the results are unchanged.) 
Thus for a matrix of order $2N+4M = 260$, we are able to
locate the desired band structure with relative error around $10^{-15}$
in the Bloch parameters $(a,b)$.
Filling such a matrix takes around 0.45 sec and computing the complex
SVD around 0.15 sec.
\footnote{All timings are reported for a laptop running \matlab\ 2008a with
  a 2GHz Intel Core Duo CPU.}
Furthermore, by storing coefficient matrices in
the expansion $\tilde{E} = \sum_{-1\le j,k \le 2}
\al^j \bt^k \tilde{E}^{(j,k)}$
at fixed $\om$, we can fill $\tilde{E}$ for new $a,b$ values in
0.05 sec.

Fig.~\ref{f:conv}b shows that, with $M$ in the new quasi-periodizing scheme
sufficient to yield machine precision,
the error convergence rate with respect to $N$ is the same
as that of the old scheme. 
Fig.~\ref{f:conv}c demonstrates the robustness of the scheme, by plotting
the smallest singular value over the same region of parameter space 
as Fig.~\ref{f:break}b.
Notice that the location of the desired band structure (black line)
is unchanged, but that the divergent behavior 
near the empty resonant band structure has entirely vanished.


\bfi 
\mbox{%
\raisebox{-2.9in}{\hspace{-2ex}\ig{height=2.9in}{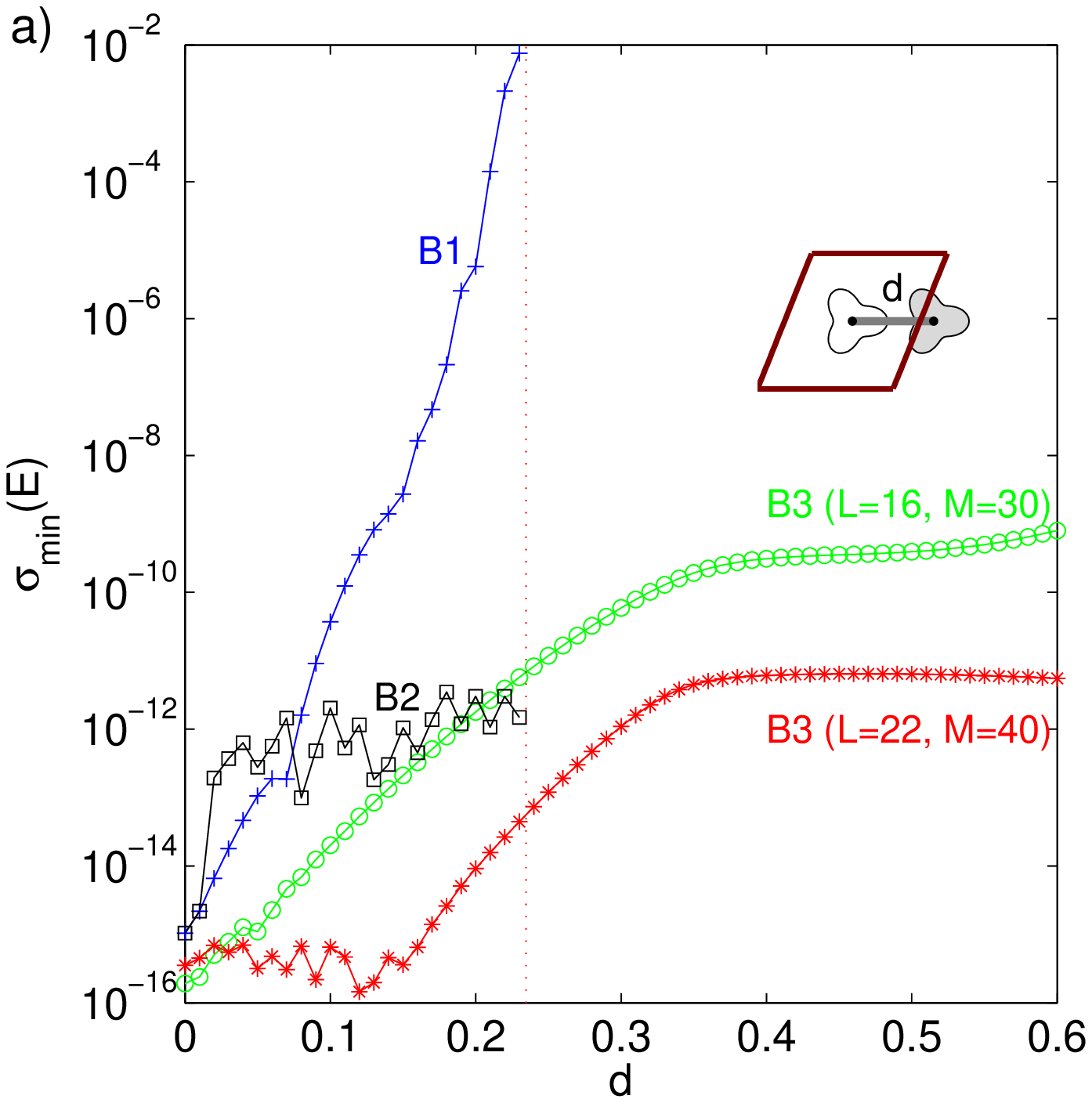}}
\;
\raisebox{-2.8in}{\hspace{-2ex}\ig{height=2.8in}%
{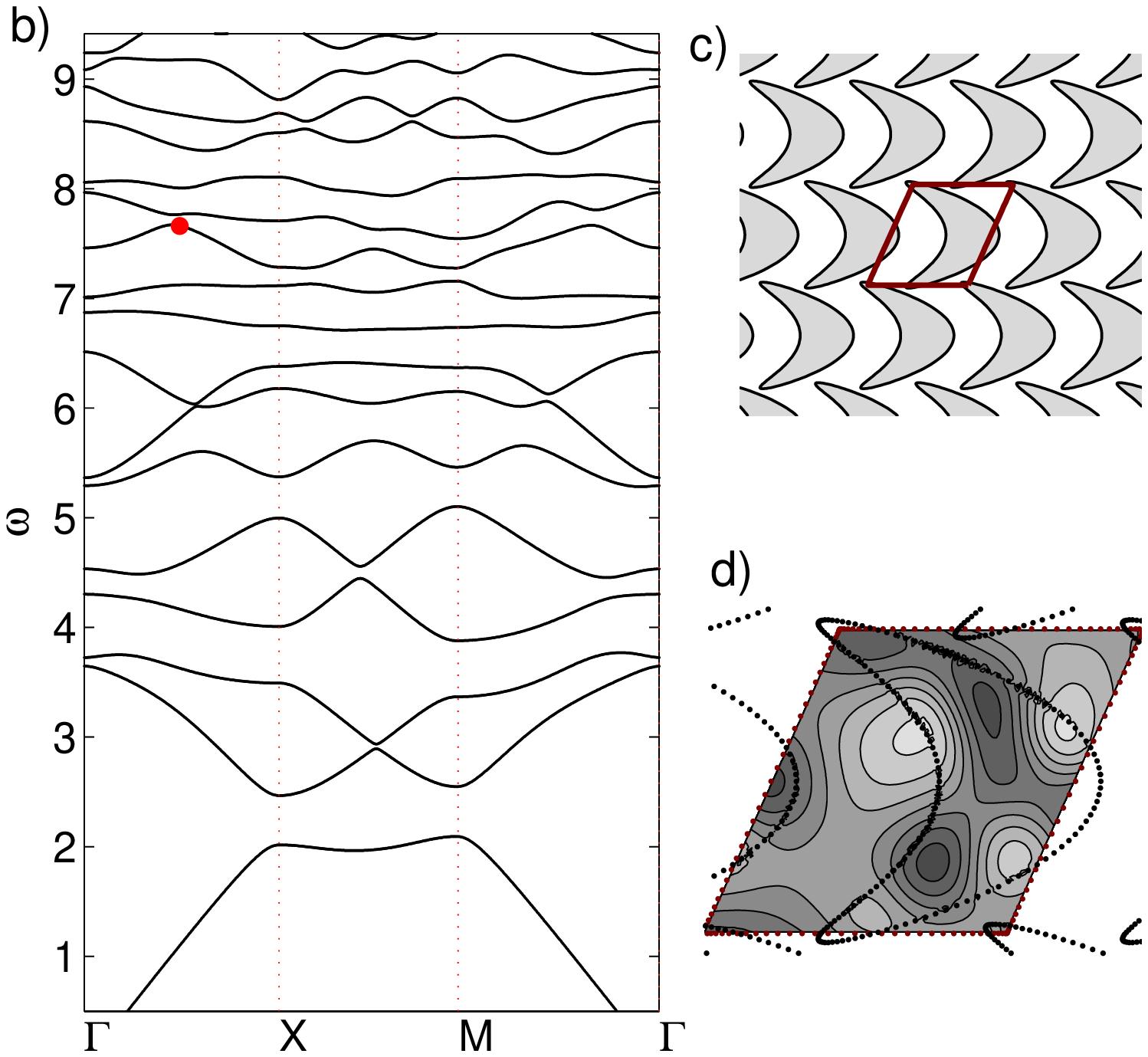}}
}%
\ca{a) Dependence of $\sigma_\tbox{min}(\tilde{E})$
on $x$-translation distance $d$ of $\pO$ relative to the system of
Fig.~\ref{f:gqpconv}a, for fixed $N=70$, and $M=30$.
The vertical line shows where $\pO$ starts to touch $\pU$.
$B$ is discretized as follows:
method B1 (+ symbols), method B2 with $d_0=0.2$ ($\square$ symbols), method B3
with $L=16$ ($\circ$ symbols), method B3 with $L=22$ and $M=40$
($\ast$ symbols).
Inset shows unit cell and inclusion at $d=0.6$.
b)
Band structure for crescent-shaped photonic crystal shown in c), index $n=2$,
shape $(0.265\cos 2\pi t +0.318\cos 4\pi t, 0.53\sin 2\pi t)$,
$0\le t<1$, unit cell $\eb=(1,0)$, $\el=(0.45,1)$.
A tour $\Gamma X M \Gamma$ of the Brillouin zone is shown, where
$\Gamma$ is $(a,b)=(0,0)$, $X$ is $(\pi,0)$, and $M$ is $(\pi,\pi)$.
d) Contours of the Bloch mode Re$[u]$ with parameters shown
by the dot on the band diagram.
}{f:wrap}
\efi

\subsection{Inclusions approaching and intersecting the unit cell wall}
\label{s:wrap}

Given a crystal of inclusions,
it may be impossible to choose a parallelogram unit cell $U$ whose boundary
does not come close to or even intersect $\pO$.
Although this is not an issue for the scheme of Section~\ref{s:gqp},
for the new scheme which relies on $\pU$ it is a potential problem.

We first show that, as expected, with method B1
the error performance deteriorates exponentially as $\pO$ approaches $\pU$.
In Fig.~\ref{f:wrap}a we plot the minimum singular value at a Bloch
eigenvalue, as a function of distance $d$ that the inclusion has been translated
in the $x$ direction (translation does not affect the Bloch eigenvalue.)
Numerical parameters $N$ and $M$ are held fixed.
The logarithm of the error grows roughly linearly
with $d$ and reaches $O(1)$ for $\dist(\pO,\pU)=0$, indicated by the dotted
vertical line at around $d=0.23$.
Method B2, also shown in Fig.~\ref{f:wrap}a, uses adaptive quadrature
for accurate evaluation of $u_\tbox{QP}$ in all of $U$.
For very small $d$, the inclusion is still centrally located (far from the wall)
and B2 is identical to B1, with an error of $10^{-15}$. The error is around
$10^{-12}$ as one approaches the wall (more or less independent of $d$),
limited by the accuracy of {\tt quadgk}.
This proves that the deterioration seen with B1 is associated with
the $B$ operator block, and can be remedied merely by careful
discretization of $B$ without increasing the matrix size.
We did not bother continuing the computation with B1 or B2 after the inclusion
crosses the wall; here they fail
because  \eqref{e:uqp}, as constructed, represents $u_\tbox{QP}$ only
inside $U$ (jump relations cause the values outside $U$ to be different).
We note that to use
B1 or B2 correctly, one would have to wrap the boundary
points outside $U$ back into the cell,
evaluate at the wrapped point, and correct for phase. 
Method B2 is not very useful in practice since the call to a black box adaptive quadrature
routine causes the matrix fill time to increase to 55 sec.

Finally we use method B3 with $L=16, M=30$ and with
$L=22, M=40$. 
In the first case, errors grow slowly to around $10^{-12}$ as $\dist(\pO,\pU)$ reaches zero,
and then continue to grow slowly to a plateau at around $10^{-9}$,
even though most of $\pO$ now falls {\em outside} of the unit cell.
The cost of B3 is not much more than B1, taking 0.7 sec to fill $\tilde{E}$.
Note that the $J$-expansion used to represent 
$u_\tbox{QP}$ has effectively carried out analytic continuation beyond $U$.
This is stable because our image structure has pushed the singularities out beyond the
nearest image cells. It is perhaps worth observing that some care must be taken in 
setting $L$. With $M=30$, increasing  $L$ above 16 would worsen errors (not shown).
The reason is that the coefficients $|l|>16$ involve more oscillatory integrands which are
not resolved by $M=30$ points. Increasing $M$ to 40 permits increased precision
with $L=22$, as seen in Fig.~\ref{f:wrap}a.

There is another potential pitfall with method B3 as implemented;
if both $L$ and $d$ get larger, there may
arise singular values of $\tilde{E}$ which become exponentially small,
associated with highly oscillatory non-physical densities on the farthest
part of $\pO$.
For illustration, with $L=16$ and $d=0.6$, the second-smallest singular value
is $10^{-4}$; with $L=22$ the second smallest singular value shrinks to  $10^{-6}$.
(When $d=0$, the second smallest singular value is  $10^{-1}$.)
This is troublesome for eigenvalue search methods that track
$\sigma_\tbox{min}(\tilde{E})$ vs Bloch parameters, since the desired minima
will be obscured by these spurious small singular values everywhere except in
a small neighborhood of the desired band structure.
We will discuss search methods less sensitive to this problem in a
future paper. For now the lesson is that, when parts of an inclusion extend
far beyond $U$, there is a price to pay for making use of analytic continuation.


\subsection{Application to band structure}

We compute the band structure of a more difficult crystal in
Fig.~\ref{f:wrap}b. $\Omega$ is far from circular,
hence simple multipole methods \cite{hafner90} would not be accurate.
The closest approach to its neighbors is only 0.06,
so that $N=150$ points are needed in discretizing the inclusion boundary.
Note that any parallelogram unit cell must intersect $\pO$,
so the method of \cite{yuan08} cannot be used without modification.
We use method B3 with $M=35$ and $L=18$.
As illustrated before in Fig.~\ref{f:gqpconv}b, the minimum values of
$\sigma_\tbox{min}(\tilde{E})$ on the band structure indicate
the size of the errors in the Bloch parameters found.
By this measure, sampling 100 random points on the first 15 bands,
we find a median error of $3\times 10^{-10}$ and a maximum $1.6\times 10^{-9}$.
1.7 sec were required to fill the matrix $\tilde{E}$ of order
$440$ once for a given
$\omega,a,b$ (and 0.13 sec for subsequent values of $a,b$).
The SVD required 0.7 sec for a matrix of this size.
We located the band structure using $8000$ such evaluations and a specialized
search algorithm, which we will describe in a forthcoming paper.
The search algorithm is also accelerated by computing the determinant
of $\tilde{E}$ rather than the SVD, at a cost of 0.1 sec for each matrix. 
The total CPU time required was 35 minutes.

Fig.~\ref{f:wrap}d shows a single Bloch mode
on the 11th band for this crescent-shaped crystal.
This took 16 sec to evaluate on a $100\times100$ grid over $U$ using
\eqref{e:urep}, and the $J$-expansion for \eqref{e:uqp} (with no 
fast multipole acceleration).

\section{Conclusions}
\label{s:conc}

We have presented two algorithms for locating the band structure
of a two-dimensional photonic crystal, in the $z$-invariant Maxwell setting.
The first (Section~\ref{s:gqp}) uses the quasi-periodic Green's function.
Theorem~\ref{t:aqp} guarantees the success of this method
(no spurious or missed modes) as long as the
band structure for the empty unit cell is avoided, where we
have shown that the method fails.
The second method (Section~\ref{s:qplp}) introduces a small number of
additional degrees of freedom on the walls to represent
the periodizing part of the field: numerical evidence suggests that
it is immune to breakdown for any Bloch parameters (Conjecture~\ref{c:qplp}).
The two schemes are connected by the following observation.
\begin{rmk}
Computing the Schur complement formula for the operator system \eqref{e:2x2}
recovers the quasi-periodic Green's function approach described by
\eqref{e:aqp}. In particular,
\[ \aqp = A - B Q^{-1} C. \]
The  quasi-periodic Green's function approach fails when $Q$ becomes singular and 
$\aqp$ blows up. The full system \eqref{e:2x2}, on the other hand, remains
well-behaved.
\end{rmk}

We have shown spectral convergence for both schemes, achieving
close to machine precision accuracy on simple crystals using only
a few hundred degrees of freedom, hence CPU times of less than 1 sec
for testing at a single parameter set $(\om,a,b)$.
In the new scheme we have shown (method B3) how to handle the passage
of the inclusion through the unit cell boundary, without much sacrifice
in accuracy, without much extra numerical effort,
and with no bookkeeping needed to determine which points of $\pO$ lie in $U$.
The latter is convenient for larger-scale or three-dimensional (3D)
computations if existing scattering codes are to be
used to fill the $A$ operator block.
Other ways to handle this intersection problem exist, such as
a variant of B2 which wraps points on $\pO$ back into $U$, with
which we have preliminary success.

We have not discussed the methods we use for the nonlinear eigenvalue
problem, due to space constraints.
The scheme of Yuan et al \cite{yuan08}
uses a quadratic eigenvalue problem,
and factorizes the scattering matrix of the inclusion at each $\om$,
hence may be faster than our scheme for small systems.
However, moving to large-scale systems with more than $10^4$ degrees of
freedom, such a factorization would be impractical compared to
an iterative version of our scheme. 

Some generalizations of what we present are straightforward, such as multiple
inclusions per unit cell, non-simply connected inclusions, or inclusions with corners
(using quadrature rules such as \cite{bremer,helsing_corner}).
There exist regimes, however, that would require some modification.
These include two phase dielectrics one or more
of which are connected through the bulk 
(sometimes called bicontinuous), and unit cells which are highly skew or have
large aspect ratios. 

Our new representation for quasi-periodic fields can also be used for 
scattering calculations from periodic one-dimensional arrays of inclusions in 2D
and one or two-dimensional arrays in 3D. Because we rely
entirely on the free-space Green's function, it should be straightforward
to create quasi-periodic solvers from existing scattering codes. 
We will describe such solvers at a later date.

\section*{Acknowledgements}
We thank Greg Beylkin, Zydrunas Gimbutas and Ivan Graham for
insightful discussions. The work of AHB was supported by 
NSF grant DMS-0811005, and by the Class of 1962 Fellowship at Dartmouth College.
The work of LG was supported by the Department of Energy under contract
DEFG0288ER25053 and by AFOSR under MURI grant FA9550-06-1-0337.

\appendix 
\section{Proof of Theorem~\ref{t:aqp}}
\label{a:aqp}

Recall the Green's representation formulae \cite[Sec. 3.2]{CK83}.
If $u$ satisfies $(\Delta + \om^2)u=0$ in $\Omega$,
recalling  that $u^-$ and $u^-_n$ signify limits on $\pO$ approaching
from the inside, and the normal always points outwards from $\Omega$, then
\be
-{\cal S}^{(\om)} u^-_n + {\cal D}^{(\om)}u^-
\; = \; \left\{\begin{array}{ll}-u& \mbox{in }\Omega\\
0 & \mbox{in } \RR\setminus\overline{\Omega}\end{array}\right.
\label{e:grfi}
\ee
The exterior representation has the opposite sign:
let $u$ satisfy $(\Delta + \om^2)u=0$ in $\RR\setminus\overline{\Omega}$
and the Sommerfeld radiation condition, that is,
\be
\frac{\partial u}{\partial r} - i\om u = o(r^{-1/2}),
\qquad  r:=|\xx| \to \infty
\label{e:rad}
\ee
holds uniformly with respect to direction $\xx/r$. Then,
\be
-{\cal S}^{(\om)} u^+_n + {\cal D}^{(\om)}u^+
\; = \; \left\{\begin{array}{ll}0& \mbox{in }\Omega\\
u & \mbox{in } \RR\setminus\overline{\Omega}\end{array}\right.
\label{e:grfe}
\ee
We will need the following quasi-periodic analogues.
\begin{lem} 
Let $u$ satisfy $(\Delta + \om^2)u=0$ in $\Omega$,
and $\overline{\Omega}\subset U$, Then for each Bloch phase $(\al,\bt)$,
\be
-{\cal S}^{(\om)}_\tbox{\rm QP} u^-_n + {\cal D}^{(\om)}_\tbox{\rm QP}u^-
\; = \; \left\{\begin{array}{ll}-u&\mbox{\rm in }\Omega\\
0 & \mbox{\rm in } U\setminus\overline{\Omega}\end{array}\right.
\ee
\label{l:grfqpi}
\end{lem}
\bp
Write $\gqp$ using \eqref{e:gqp} and
notice that each term other than
$(m,n)=(0,0)$ contributes zero. This is because all points in
$U$ lie outside each closed curve $\pO - m\eb- n\el$,
and we may apply the second (extinction) case of \eqref{e:grfi} to show
that they have no effect in $U$.
\ep
\begin{lem} 
Let $u$ satisfy $(\Delta + \om^2)u=0$ in $U\setminus \overline{\Omega}$
and quasi-periodicity \eqref{e:f}-\eqref{e:gp}, and
$\overline{\Omega}\subset U$. Then
\be
-{\cal S}^{(\om)}_\tbox{\rm QP} u^+_n +  {\cal D}^{(\om)}_\tbox{\rm QP}u^+
\; = \;  \left\{\begin{array}{ll}0&\mbox{\rm in }\Omega\\
u & \mbox{\rm in } U\setminus\overline{\Omega}\end{array}\right.
\ee
\label{l:grfqpe}
\end{lem}
\bp We follow the usual method of proof \cite[Thm. 3.3]{CK83} but with the
quasi-periodicity condition playing the role of the
radiation condition. Apply Green's 2nd identity to the
functions $u$ and $\gqp(\xx,\cdot)$ in the domain
$U\setminus \overline{\Omega}$ if $\xx\in\Omega$,
or the domain $\{\yy\in U\setminus \overline{\Omega}: |\xx-\yy|>\eps\}$
if $\xx\in U\setminus \overline{\Omega}$.
In the latter case the limit $\eps\to0$ is taken,
and \eqref{e:gqp} shows that only the $(m,n)=(0,0)$ term contributes to the
limit of the integral over the sphere of radius $\eps$.
In both cases the boundary integrals contain the term
\be
\int_{\partial U} \frac{\partial \gqp}{\partial n_\yy}(\xx,\yy)u(\yy)
- \gqp(\xx,\yy)u_n(\yy) \;ds_\yy~,
\ee
which vanishes by cancellation on opposing walls, since
$u$ is quasi-periodic with phases $(\al,\bt)$, but 
$\gqp(\xx,\cdot)$ is anti-quasiperiodic, i.e.\ quasi-periodic
with phases $(\al^{-1},\bt^{-1})$.
\ep

Turning now to Theorem~\ref{t:aqp},
to prove the {\em if} part, we show that whenever the operator has a
nontrivial nullspace, a Bloch eigenfunction $u$ may be constructed,
i.e.\ a solution to \eqref{e:pdei}-\eqref{e:gp} that we must take care to
show is nontrivial.
Let $\eta = [\tau;-\sigma] \neq 0$ be a nontrivial density such that
$\aqp \eta = 0$.
Immediately we have that the resulting field $u$ given by \eqref{e:qprep}
satisfies  \eqref{e:pdei}-\eqref{e:gp}.
We now define a complementary field over the whole plane minus $\pO$,
\be
v \;=\; \left\{\begin{array}{ll}
{\cal S}^{(\om)}_\tbox{QP}\sigma + {\cal D}^{(\om)}_\tbox{QP}\tau &
\mbox{in } \Omega\\
-{\cal S}^{(n\om)}\sigma - {\cal D}^{(n\om)}\tau & \mbox{in }
\RR\setminus\overline{\Omega}
\end{array}\right.
\label{e:v}
\ee
Suppose $u\equiv 0$. Then $u^- = u_n^- = 0$ and by the jump relations for
${\cal S}^{(n\om)}\sigma + {\cal D}^{(n\om)}\tau$ we get
$v^+ = -\tau$ and $v_n^+ = \sigma$.
Similarly, since  $u^+ = u_n^+ = 0$ by the jump relations for
${\cal S}^{(\om)}_\tbox{QP}\sigma + {\cal D}^{(\om)}_\tbox{QP}\tau$
we get $v^- = -\tau$ and $v_n^- = \sigma$.
It is easy to check that $v$ solves the
(swapped-wavenumber) transmission problem,
\bea
(\Delta+\omega^2)v &=& 0 \qquad \mbox{in } \Omega 
\label{e:trani}
\\
(\Delta+n^2\omega^2)v &=& 0 \qquad \mbox{in } \RR \setminus \overline{\Omega}
\label{e:trane}
\\
\frac{\partial v}{\partial r} - i n \om v &=& o(r^{-1/2}),
\qquad  r \to \infty, \mbox{ uniformly in direction}
\\
v^+-v^-&=& h
\\
v_n^+-v_n^-&=& h'
\label{e:Hp}
\eea
with homogeneous boundary discontinuity data $h=h'=0$.
By uniqueness for this problem
\cite[Thm.~3.40]{CK83} we get that $v\equiv 0$ in $\RR$, from which
the jump relations back to $u$ imply $\sigma = \tau = 0$,
which contradicts our assumption of nontrivial density.
Thus $u$ is a Bloch eigenfunction.

To prove the {\em only if} part
we show that, given the existence of a Bloch eigenfunction,
we may exhibit a (nontrivial) density $\eta$ such that $\aqp \eta = 0$.
Let $w$ be a Bloch eigenfunction with eigenvalue $(\om,a,b)$.
Then let $v$ solve \eqref{e:trani}-\eqref{e:Hp} with the inhomogeneous
data $h = -2 w|_\pO$ and $h' = -2 w_n|_\pO$.
(Note that $w$ obeys continuity \eqref{e:cont},
hence $w|_\pO=w^+=w^-$ and $w_n|_\pO = w^+_n = w^-_n$).
By \cite[Thm.~3.41]{CK83} we know that a unique solution exists.
We now claim that the densities
\bea
\sigma &=& w_n|_\pO + v_n^+
\label{e:sig}
\\
\tau &=& -w|_\pO - v^+
\label{e:tau}
\eea
generate precisely the eigenfunction $w$, i.e.\
the representation $u$ of \eqref{e:qprep} obeys $u\equiv w$ in $U$.
We show this by substituting the densities into \eqref{e:qprep},
then applying \eqref{e:grfi} and \eqref{e:grfe} in $\Omega$,
and Lemma~\ref{l:grfqpe} in $U\setminus\overline{\Omega}$:
\bea
u &=& \left\{\begin{array}{ll}
{\cal S}^{(n\om)}w_n|_\pO - {\cal D}^{(n\om)}w|_\pO
\;\;+ \;\;
{\cal S}^{(n\om)}v_n^+ - {\cal D}^{(n\om)}v^+ & \mbox{in } \Omega\\
{\cal S}^{(\om)}_\tbox{QP}w_n|_\pO - {\cal D}^{(\om)}_\tbox{QP}w|_\pO
\;\;+\;\;
{\cal S}^{(\om)}_\tbox{QP}v_n^+ - {\cal D}^{(\om)}_\tbox{QP}v^+ & \mbox{in }
U\setminus\overline{\Omega}
\end{array}\right.
\nonumber\\ 
&=&
 \left\{\begin{array}{ll}
w & \mbox{in } \Omega\\
-w \;\;+\;\;
{\cal S}^{(\om)}_\tbox{QP}v_n^+ - {\cal D}^{(\om)}_\tbox{QP}v^+
&\mbox{in } U\setminus\overline{\Omega}
\end{array}\right.
\nonumber
\eea
On the remaining term,
we use $v$'s known jumps $h$ and $h'$ to get
\bea
{\cal S}^{(\om)}_\tbox{QP}v_n^+ - {\cal D}^{(\om)}_\tbox{QP}v^+
&=&
{\cal S}^{(\om)}_\tbox{QP}v_n^- - {\cal D}^{(\om)}_\tbox{QP}v^-
- 2 {\cal S}^{(\om)}_\tbox{QP}w_n|_\pO + 2 {\cal D}^{(\om)}_\tbox{QP}w|_\pO
\nonumber\\
&=& -2w \nonumber
\eea
where we applied Lemma~\ref{l:grfqpi} to the first pair, and
Lemma~\ref{l:grfqpe} to the second as before.
Substituting this above shows that $u\equiv w$ in $U$.
Since $w$ has zero mismatch, the density vector $\eta := [\tau;-\sigma]$
satisfies $\aqp \eta = 0$.
Finally, $\eta$ must be nontrivial since $\eta=0$ would imply $u\equiv 0$
by \eqref{e:qprep} which contradicts it being equal to the eigenfunction $w$.
\ep

We close with a couple of remarks about the proof.
Barring a sign, $v$ in \eqref{e:v} is the extension of $u$'s representation
\eqref{e:qprep} into its nonphysical regions, a 
trick originating, in the homogeneous context,
with the proof in \cite[Thm.~3.41]{CK83}.
Because \eqref{e:qprep} uses $\gqp$ outside, but $G$ inside,
the complementary problem is a {\em nonperiodic}
transmission problem, which has known existence and uniqueness.
The related analysis of \cite{shipman} uses $\gqp$ both inside
and outside. This results in a periodic problem as the complementary
problem, and it is not so clear that one can eliminate the possibility of spurious modes.



\bibliography{alex}

\end{document}